\documentclass[draft]{tran-l}

\usepackage{amsfonts,amssymb}
\usepackage{graphicx}
\usepackage{checkend}
\usepackage{longtable}
\usepackage{physics}
\usepackage{booktabs} 
\usepackage{cleveref}
\usepackage{etoolbox} 
\usepackage{subfig}
\usepackage{xcolor}
\usepackage{mathrsfs}
\usepackage{enumerate}
\usepackage{mathtools}
\usepackage[utf8]{inputenc}
\usepackage[english]{babel}

\newcommand{\bb}[1]{\mathbb{#1}}

\newcommand{\nc}{\newcommand}
\nc{\ep}{\varepsilon}
\nc{\n}[1]{\mathscr{#1}}
\nc{\eps}[1]{{#1}_{\varepsilon}}
\nc{\be}{\begin{equation}}
\nc{\ee}{\end{equation}}
\nc{\m}[1]{\mathcal{#1}}
\DeclareMathOperator*{\esssup}{ess\,sup}

\newtheorem{theorem}{Theorem}[section]
\newtheorem{lemma}[theorem]{Lemma}

\theoremstyle{definition}
\newtheorem{definition}[theorem]{Definition}
\newtheorem{corollary}{Corollary}

\theoremstyle{remark}

\numberwithin{equation}{section}

\begin{document}

\title[Optimal Control of free boundary problems]{Optimal control of multiphase free boundary problems for nonlinear parabolic equations}


\author[U.G. Abdulla]{Ugur G. Abdulla}
\address{Department of Mathematical Sciences, Florida Institute of Technology, 150 W University Blvd, Melbourne FL, 32901}
\curraddr{}
\email{abdulla@fit.edu}
\thanks{}

\author[E. Cosgrove]{Evan Cosgrove}
\address{Department of Mathematical Sciences, Florida Institute of Technology, 150 W University Blvd, Melbourne FL, 32901}
\curraddr{}
\email{ecosgrove2011@my.fit.edu}
\thanks{}

\keywords{Inverse multidimensional multiphase Stefan problem, Quasilinear parabolic
PDE with discontinuous coefficients, optimal control, Sobolev spaces, method of finite differences, discrete optimal control problem, energy estimate, embedding theorems, weak compactness, convergence in functional, convergence in control, maximal monotone graph}
\subjclass[2010]{35R30, 35R35, 35K20, 35Q93, 49J20, 65M06, 65M12, 65M32,
65N21}

\date{}

\dedicatory{}

\begin{abstract}

  We consider the optimal control of singular nonlinear partial differential equation which is the distributional formulation of the multiphase Stefan type free boundary problem for the general second order parabolic equation.  Boundary heat flux is the control parameter, and the optimality criteria consist of the minimization of the $L_2$-norm declination of the trace of the solution to the PDE problem at the final moment from the given measurement. Sequence of finite-dimensional optimal control problems is introduced through finite differences. We establish existence of the optimal control and prove the convergence of the sequence of discrete optimal control problems to the original problem both with respect to functional and control. Proofs rely on establishing a uniform $L_{\infty}$ bound, and $W_2^{1,1}$-energy estimate for the discrete nonlinear PDE problem with discontinuous coefficient.

\end{abstract}

\maketitle
\section{Introduction}
\subsection{Optimal Control Problem}
Let $v^1<v^2<\cdots<v^m$ are given real numbers. Consider an optimal control problem on the minimization of the cost functional
\begin{equation}\label{functional}
    \n J(g) = \Vert v(x,T;g)-\omega(x)\Vert^2_{L_2(0,\ell)}
\end{equation}
over the control set:

\[\n G_R =\big\{g: g\in W_2^1(0,T),\,\,\Vert g\Vert_{W_2^1[0,T]}\leq R\big\}.\]
where $v$ is a solution of the singular nonlinear PDE problem
\begin{align}
    \frac{\partial \beta(v)}{\partial t}-\mathcal{L}v-f(x,t) \ni 0, &\qquad 0<x<\ell,\ 0<t\leq T; \ v(x,t)\neq v^j, j=\overline{1,m} \label{bvpde}\\
    v(x,0)=\Phi(x) \label{vphi}, &\qquad 0 < x < \ell\\
    av_x+bv|_{x=0}=g(t),&\qquad 0<t\leq T,  \label{vg}\\
    av_x+bv|_{x=\ell} = p(t),&\qquad 0<t\leq T,  \label{vp}
\end{align}
with $\beta(\cdot)$ being a maximal monotone graph of the form
\begin{equation}\nonumber
    \beta(r)=
    \left\{
    \begin{array}{l}
    \beta_j(r)+\sum\limits_{i=0}^{j-1}\nu_i, \quad\text{for} \ v^{j-1}<r<v^j,\\
    $\Big[$\beta_j(v^j)+\sum\limits_{i=0}^{j-1}\nu_i,\beta_j(v^j)+\sum\limits_{i=0}^{j}\nu_i $\Big]$, \quad\text{for} \  \ r=v^j,\\
    \beta_{j+1}(r)+\sum\limits_{i=0}^{j}\nu_i, \quad\text{for} \ v^j<r<v^{j+1}; \ j=1,2,...,m
    \end{array}\right.
\end{equation}
with a given positive constants $\nu_j, j=1,...,m$; $\nu_0=0$, $v^0=-\infty$, $v^{m+1}=+\infty$;  $\beta_i(\cdot), i=1,...,m+1$ are monotone increasing Lipschitzian functions in their respective domain of definition, $\beta_j(v^j)=\beta_{j+1}(v^j), j=1,...,m$, 
    \begin{equation}\label{bbar} 
        0<\bar{b} \leq \beta_j^{\prime}(r) \leq \bar{B}, \ j=1,...,m+1;
    \end{equation}
and $\mathcal{L}$ is an elliptic operator
\[ \mathcal{L}v=(a(x,t)v_x+b(x,t)v)_x-c(x,t)v \]
with bounded measurable coefficients $a,b,c$ and
\begin{equation}\label{apositive}
a(x,t)\geq a_0>0, \ \quad\text{a.e.} \ (x,t)\in D=\{0<x<\ell, 0<t<T\}.
\end{equation}
Described optimal control problem will be called a Problem $\mathcal{S}$. 

In the particular case of $\mathcal{L}=\Delta$, system \eqref{bvpde}-\eqref{vp} presents distributional formulation of the multiphase Stefan problem describing flow of the heat in the presence of phase transitions \cite{Oleinik, Kamenomostskaya2, LSU}. In the physical context, $v(x,t)$ is a temperature, $f(x,t)$ is a density of heat sources, $\Phi(x)$ is an initial temperature distribution, $g(t)$ and $p(t)$ are heat flux on fixed boundaries, $v^j$'s are phase transition temperatures; $\beta_j^{\prime}(v), v^j<v<v^{j+1}, j=0,1,...,m$ characterize heat conductivities in each phase, and the positive jump constants $\nu^j, j=1,...,m$ are expressing latent heat of fusion during phase transition. In particular, the case $m=1, v^1=0$ is a classical two-phase Stefan problem describing melting of the ice or freezing of the water. \cite{Friedman1,Meyrmanov}. More complex examples of multiphase Stefan problem includes biomedical problem about the laser ablation of biomedical tissues, which motivates general elliptic operator $\mathcal{L}$ with coefficients $a,b,c$ expressing anisotropic properties of the media. Optimal control problem \eqref{functional}-\eqref{vp} aims to achieve the desired temperature distribution $\omega(x)$ at the final moment by controlling the boundary flux $g(t)$ on the fixed boundary. Equivalently, it is a variational formulation of the inverse multiphase Stefan type free boundary problem on the identification boundary flux $g(t)$ through measurement $\omega(x)$ of the final moment temperature distribution. 

The goal of this paper is to prove the well-posedness of the optimal control problem \eqref{functional}-\eqref{vp}, and to prove the convergence of the sequence of the finite-dimensional discretized optimal control problems to the optimal control problem \eqref{functional}-\eqref{vp} both with respect to functional and control via the method of finite differences. 
    
The idea of transformation of the multiphase Stefan problem to boundary value problem for singular PDE \eqref{bvpde} originated in \cite{Oleinik}. 
Existence and uniqueness of the weak solution was proved in \cite{Oleinik, Kamenomostskaya2,LSU} when $\mathcal{L}=\Delta$. H\"{o}lder continuity of the weak solutions 
was proved in \cite{DiBenedetto2,DiBenedetto} for general nonlinear elliptic operators $\mathcal{L}$. Continuity of the weak solution for the two-phase Stefan problem was proved in \cite{Caffarelli}.

The one-phase inverse Stefan problem (ISP) was first mentioned in \cite{Cannon3}, where phase transition boundary is known and heat flux on left boundary is to be found, and the variational approach for solving the ISP was used in \cite{BudakVasileva1, BudakVasileva2}. In \cite{Vasilev} ISP was formulated as an optimal control problem and the existence of the optimal control is proved. In \cite{Yurii}, the Frechet derivative was found, the convergence of finite difference schemes was proved, and Tikhonov regularization was suggested in order to improve results. The following works on the ISP split into two different directions: ISPs with given fixed phase transition boundaries (\cite{Alifanov,Bell,Budak,Cannon,Carasso,Ewing1,Ewing2,Hoffman,Sherman,Goldman}), and ISPs with unknown phase transition boundaries (\cite{Baumeister,Fasano,Hoffman1,Hoffman2,Jochum2,Jochum1,Knabner,Lurye,Nochetto, Niezgodka,Primicero,Sagues,Talenti,Goldman}). We can refer to \cite{Goldman} for a full list of references for both types stated above, which include both linear and quasilinear parabolic equations.

In \cite{Abdulla1, Abdulla2} a new variational formulation of the the one-phase ISP was developed, in which optimal control framework was implemented, where the phase transition boundary is included in the control set along with the boundary heat flux. The sum of the $L_2$-norm declinations are minimized against the available measurements of temperature on the fixed boundary, available measurements of the free boundary location, and temperature at final moment.  Important advantage of the new control theoretic approach is that it can handle situations where the phase transition temperature is not known explicitly, and is available only through measurement with possible error. Another major advantage of the new variational method suggested in \cite{Abdulla1, Abdulla2} is based on the fact that for a given control vector corresponding state vector solves PDE problem in a fixed region instead of full free boundary problem. This allows to reduce significantly computational cost of iterative numerical methods based on gradient type methods in Sobolev spaces. In \cite{Abdulla3,Abdulla4}, Frechet differentiability in Sobolev-Besov spaces was proved and the formula for the Frechet gradient and optimality condition are derived. In \cite{Abdulla6,Abdulla8} gradient method was implemented in Hilbert-Besov spaces framework for the numerical solution of the ISP.

The new method developed in \cite{Abdulla1,Abdulla2} is not applicable to inverse multiphase free boundary problems. The reason is that the Stefan condition on the free boundary includes the saltus of the boundary flux from neighbouring phases, and   by fixing free boundary as a control parameter the Stefan condition does not become a Neumann or Robin type boundary condition for the PDE. In a recent paper \cite{Abdulla5}, a new variational method was introduced for the solution of the inverse multiphase Stefan problem. The IMSP is reformulated in a new optimal control framework  in which the boundary is fixed, yet the state vector satisfies a nonlinear PDE with coefficients possessing jump discontinuities along phase transition boundaries. In \cite{Abdulla5} existence of the optimal control and convergence of the sequence of discretized optimal control problems via method of finite differences is proved. In \cite{Abdulla7}, this framework is extended to the multidimensional IMSP.

\subsection{Notation of Function Spaces}

$L_2(0,T)$ - Space of Lebesgue square-integrable functions. It is a Hilbert space with inner product
\[ (u,v)=\int_0^T uv\,dt. \]
$L_2(D)$ - Hilbert space with inner product
\[ (u,v)=\int_{D} uv\, dx\, dt . \]
$L_{\infty}(0,T)$ - Space of essentially bounded functions. It is a Banach space with norm
\[
\Vert u\Vert_{L_{\infty}[0,T]} =\underset{0\leq t\leq T}{\text{ esssup }}|u(t)|.
\]
$L_{\infty}(D)$ - Space of essentially bounded functions. It is a Banach space with norm
\[
\Vert u\Vert_{L_{\infty}(D)} =\underset{(x,t)\in D}{\text{ esssup }}|u(x,t)|.
\]
$W_2^k(0,T), k=1,2,...$ - Hilbert space of all elements of $L_2(0,T)$ whose weak derivatives up to order $k$  exist and belong to $L_2(0,T)$. The inner product is defined as
\[ (u,v)=\int_0^T \sum_{s=0}^k \frac{d^su}{dt^s} \frac{d^sv}{dt^s} \,dt.  \]
$W_2^{1,0}(D)$ - Hilbert space of all elements of $L_2(D)$ that have a weak derivative in the $x$ direction, $\frac{\partial u}{\partial x}$, and such that it belongs to $L_2(D)$. The inner product is defined as
\[ (u,v)=\int_{D} \Big ( uv + \frac{\partial u}{\partial x}\frac{\partial v}{\partial x} \Big ) \,dx\, dt. \]
$W_2^{1,1}(D)$ - Hilbert space of all elements of $L_2(D)$ with weak derivatives of first order, $\frac{\partial u}{\partial x}$, $\frac{\partial u}{\partial t}$. Also its weak derivatives must belong to $L_2(D)$. The inner product is defined as
\[ (u,v)=\int_{D} \Big ( uv + \frac{\partial u}{\partial x}\frac{\partial v}{\partial x}+ \frac{\partial u}{\partial t}\frac{\partial v}{\partial t} \Big )\, dx\,dt. \]
$W_\infty^{1,0}(D)$ - Space of elements of $L_\infty(D)$ with weak derivative in the $x$ direction existing and belonging to $L_\infty(D)$. It is a Banach space with norm
\[\|u\|_{W_\infty^{1,0}(D)} = \|u\|_{L_\infty(D)} + \|u_x\|_{L_\infty(D)} \]

\subsection{Weak Solution of the Multiphase Free Boundary Problem}

We now formulate the notion of the weak solution of the nonlinear multiphase parabolic free boundary problem (\ref{bvpde})-(\ref{vp}).

\begin{definition}\label{typeB} We say that a measurable function $B(x,t,v)$ is \textit{of type }$\n B$ if
\begin{enumerate}[(a)]
	\item $B(x,t,v)=\beta(v),\qquad v\neq v^j,\quad\forall j=\overline{1,J}$
	\item $B(x,t,v)\in[\beta(v^j)^-,\beta(v^j)^+],\qquad v=v^j$ for some $j$. \\
\end{enumerate}
\end{definition}
\begin{definition}\label{weaksol} $v\in W_2^{1,1}(D)\cap L_{\infty}(D)$ is called a \emph{weak solution of the problem}  (\ref{bvpde})-(\ref{vp}) if for any two functions $B,B_0$ of type $\n B$, the following integral identity is satisfied
\begin{gather}
\int\limits_0^T\int\limits_0^{\ell}\Big[-B(x,t,v(x,t))\psi_t+a(x,t)v_x\psi_x+b(x,t)v\psi_x+c(x,t)v\psi-f\psi\Big]\,dxdt \nonumber
\\
- \int\limits_0^{\ell}B_0(x,0,\Phi(x))\psi(x,0)\,dx
- \int\limits_0^Tp(t)\psi(\ell,t)\,dt
+\int\limits_0^Tg(t)\psi(0,t)\,dt=0,\label{weaksolution}
\end{gather}
for all \(\psi\in W_2^{1,1}(D)\) with \(\psi(x,T)=0\).
\end{definition}

\subsection{Discrete Optimal Control Problem}

Let 
\[\omega_{\tau} = \{t_k,k=\overline{1,n}\},\tau=\frac Tn, t_k=k\tau,\qquad\quad \omega_h= \{x_i,i=\overline{1,m}\},h=\frac{\ell}m, x_i=ih
\]
be grids in the time and space domains, respectively, under the assumptions that
$m\rightarrow\infty\quad$ as $n\rightarrow\infty$ and
\begin{gather}\label{htau}
\frac{h}{\tau} \geq \frac{8\|b\|_{L_\infty}}{\bar{b}}
\end{gather}
Define the Steklov averages
\begin{gather}\label{steklov}
w_k = \frac1{\tau}\int\limits_{t_{k-1}}^{t_k}w(t)\,dt,\qquad\qquad \Phi_i = \frac 1h\int\limits_{x_i}^{x_{i+1}}\Phi(x)\,dx,~~ \Phi_m=\Phi(\ell), \\ q_{ik}=\frac{1}{\tau h}\int\limits_{t_{k-1}}^{t_k} \int\limits_{x_i}^{x_{i+1}}q(x,t)\,dxdt,\qquad k=\overline{1,n},\quad i=\overline{0,m-1} \nonumber,
\end{gather}
where $w$ represents any of the functions $p$, $\Gamma$, $g$, or $g^n$, and $q$ represents any of the functions $a, b, c,$ and $f$. Introduce the discretized control set
\[\n G_R^n = \{[g]_n\in\bb R^{n+1}: \Vert[g]_n\Vert_{w_2^1}\leq R\}\]
where $[g]_n=(g_0,g_1,\ldots,g_n)$, and
\[
\Vert[g]_n\Vert_{w_2^1}^2=\sum\limits_{k=1}^{n}\tau g_k^2+\sum\limits_{k=1}^n\tau g_{k\bar t}^2
\]
with $g_{k\bar t}=\frac{g_k-g_{k-1}}{\tau}.$ Assume that every element $g\in W_2^1(0,T)$ is extended on the interval $[-\tau, 0]$ as a constant $g(0)$. Consider now the mappings between the discrete and continuous control sets, $\n Q_n: W_2^1(0,T)\rightarrow \bb R^{n+1},\quad \n P_n:\bb R^{n+1}\rightarrow W_2^1(0,T)$ as
\begin{gather}
\n Q_n(g)=[g]_n,\quad\text{for }g\in\n G_R,\,\text{where }g_k=\frac1{\tau}\int\limits_{t_{k-1}}^{t_k}g(t)\,dt,\quad k=\overline{0,n},\\
\n P_n([g]_n)=g^n,~\text{for }[g]_n\in\n G_R^n;\nonumber \\
g^n(t)=g_{k-1}+\frac{g_k-g_{k-1}}{\tau}(t-t_{k-1}), t\in[t_{k-1},t_k),\,\, k=\overline{1,n}. \label{Pmap}
\end{gather}
Approximate the function $\beta(v)$ by the infinitely differentiable sequence 
\begin{equation}\label{smoothing}
b_{n}(v)=\int_{v-\frac{1}{n}}^{v+\frac{1}{n}}\beta(y)\omega_n(v-y)dy,
\end{equation}
where $\omega_n$ is a standard mollifier defined as
\begin{equation}\label{kernel}
\omega_n(v) =\left\{\begin{matrix}\m C n e^{-\frac{1} {1-n^2v^2}},\quad&|v|\leq\frac{1}{n}\\[2mm] 0,\quad&|v|>\frac{1}{n}\end{matrix}\right.
\end{equation}
and the constant $\m C$ is chosen so that $\int\limits_{\bb R}\omega_1(u)\,du =1$. Since $\beta'(v)$ is piecewise-continuous, we also have
\begin{equation}\label{smoothing_derivative}
b_{n}'(v)=\int_{v-\frac{1}{n}}^{v+\frac{1}{n}}\beta'(y)\omega_n(v-y)dy.
\end{equation}
This implies $b_n$ is also strict monotonically increasing and by \eqref{bbar} we have
\begin{equation}\label{b_n bound}
b_n'(v)\geq \bar{b}>0
\end{equation}
We now define a solution to the problem \eqref{bvpde}-\eqref{vp} in the discrete sense\\

\textbf{Discrete State Vector.} Given $[g]_n$, the vector function $[v([g]_n)]_n\\=\big(v(0),v(1),\ldots,v(n)\big);$ $\quad v(k)\in\bb R^{m+1},\quad k=0,\ldots,n$ is called a \emph{discrete state vector} if
\begin{enumerate}[(a)]
	\item $v_i(0) = \Phi_i,\quad i =\overline{0,m}, $\\
	\item For arbitrary $k=1,\ldots,n$, the vector $v(k)\in\bb R^{m+1}$ satisfies
	\begin{gather}\label{dsvsum}
	\sum\limits_{i=0}^{m-1}h\Big[\big(b_n(v_i(k))\big)_{\bar t} \eta_i+a_{ik}v_{ix}(k)\eta_{ix}+b_{ik}v_{i}(k)\eta_{ix}+c_{ik}v_i(k)\eta_i-f_{ik}\eta_i\Big]
	\\
	-p_k\eta_m+g_k^n\eta_0=0,
	\forall \eta=(\eta_i)\in\bb R^{m+1}.\nonumber
	\end{gather}
\end{enumerate}
Given $[g]_n\in\n G^n_R$, the discrete cost functional $\n I_n$ is defined as
\begin{equation}
\n I_n([g]_n) = \sum\limits_{i=1}^m h\Big(v_i(n)-w_i\Big)^2 \label{In}
\end{equation}
where $v_i(k)$ are components of the discrete state vector $[v([g]_n)]_n$. Finite-dimensional optimal control problem
on the minimization of $\n I_n([g]_n)$ on a control set $\n G_R^n$ will be called Problem $\mathcal{S}_n$.  We define 
\[\n I_{n_*}:= \inf\limits_{[g]_n\in\n G_R^n}\n I_n([g]_n).\]
Furthermore, the following interpolations will be considered:
\begin{gather}
\tilde v(x,t) = v_i(k),\qquad x\in[x_i,x_{i+1}],\quad t\in[t_{k-1},t_k],\qquad i =\overline{0,m-1},\quad k =\overline{0,n}, \nonumber \\ 
\hat v(x;k) = v_i(k)+v_{ix}(k)(x-x_i),\qquad x\in[x_i,x_{i+1}],\quad i=\overline{0,m-1}, \nonumber \\
v^{\tau}(x,t) = \hat v(x;k),\qquad t\in[t_{k-1},t_k], \nonumber \\
\hat v^{\tau}(x,t) = \hat v(x;k-1) +\hat v_{\bar t}(x;k)(t-t_{k-1}),\qquad t\in[t_{k-1},t_k],\quad k = \overline{1,n}. \label{interpolations}
\end{gather}

\section{Main Results}

Unless stated otherwise, throughout the paper we assume the following conditions are satisfied by the data:

\begin{gather}
    f\in L_{\infty}(D),\quad p\in W_2^1(0,T),\quad\Phi\in W_2^1(0,\ell), \quad a,b \in W_\infty^{1,0}(D), \quad c \in L_\infty(D), \nonumber \\
    \int\limits_{0}^{T}\esssup\limits_{x\in[0,\ell]}|a_t(x,t)|\,dt < \infty \quad \int\limits_{0}^{T}\esssup\limits_{x\in[0,\ell]}|b_t(x,t)|\,dt < \infty \label{assumptions}
\end{gather}
the coefficient $a$ satisfies \eqref{apositive}; and $\Phi(x)=v^j, j=1,...,m$ on a set of measure 0 in the interval $(0,\ell)$.

\begin{theorem}\label{existence}
	The optimal control problem $\mathcal{S}$ has a solution, that is, the set 
	\[\n G_*=\Big\{g\in\n G_R\big|\n J(g)=\n J_*:=\underset{g\in \n G_R}{\inf }\n J(g)\Big\}\]
	is not empty.
\end{theorem}

\begin{theorem}\label{convergence}
	The sequence of discrete optimal control problems $\mathcal{S}_n$ approximates the optimal control problem $\mathcal{S}$ with respect to functional, that is,
	\begin{equation}\label{Eq:W:1:18}
	\lim\limits_{n\to +\infty} \n{I}_{n_*}=\n{J}_*, 
	\end{equation}
	where
	\[ \n{I}_{n_*}=\inf\limits_{\n G_R^n} \n{I}_n([g]_n), \ n=1,2,\ldots ~.\]
	If $[g]_{n_\ep}\in \n G_R^n$ is chosen such that
	\[ \n{I}_{n_*} \le \n{I}_n([g]_{n_\ep})\le \n{I}_{n_*}+\ep_n, \ \ep_n \downarrow 0, \]
	then the sequence $g^n=\n{P}_n([g]_{n_\ep})$ has a subsequence convergent to some element $g_*\in\n G_*$ weakly
	in $W_2^1(0,T)$ and strongly in $L_2(0,T)$. Moreover, the piecewise linear interpolations $\hat{v}^\tau$ of the corresponding discrete state vectors $[v([g]_{n_\ep})]_n$ converge to the weak solution $v(x,t;g_*) \in W_2^{1,1}(D)\cap L_\infty(D)$ of the singular PDE problem (\ref{bvpde})-(\ref{vp}) weakly in $W_2^{1,1}(D)$, strongly in $L_2(D)$, and almost everywhere on $D$.
\end{theorem} 

\section{Preliminary Results}\label{PR}

\begin{lemma}\label{eudsv} 
    Given any $[g]_n\in\n G^n$, and any $h,\tau$, a discrete state vector exists uniquely.
\end{lemma}

\noindent\textit{Proof.} First we prove uniqueness by induction. For a given $[g]_n$, suppose $v$ and $\tilde v$ both are discrete state vectors. Due to definition of how a discrete state vector is constructed, we have that $v(0)=\tilde v(0)$. Now suppose that $v(k-1)=\tilde v(k-1)$ for some fixed $k \geq 1$. Since $v$ and $\tilde v$ both satisfy (\ref{dsvsum}), subtract the identities for both $v$ and $\tilde{v}$, choosing $ \eta=v(k)-\tilde v(k)$ to get:
\begin{gather}
\sum\limits_{i=0}^{m-1}\Big[\big(b_n(v_i(k))_{\bar t}-b_n(\tilde v_i(k))_{\bar t}\big)\big(v_i(k)-\tilde v_i(k)\big)+a_{ik}\big(v_{ix}(k) -\tilde v_{ix}(k)\big)^2 \nonumber
\\
+b_{ik}\big(v_{ix}(k) -\tilde v_{ix}(k)\big)\big(v_{i}(k) -\tilde v_{i}(k)\big)
+c_{ik}\big(v_{i}(k) -\tilde v_{i}(k)\big)^2 \Big]=0.
\end{gather}

From here we get using Cauchy inequality with $\epsilon = a_0$:

\begin{gather}
\sum\limits_{i=0}^{m-1}\Big[\frac{1}{\tau}\big(b_n(v_i(k))-b_n(\tilde v_i(k))\big)\big(v_i(k)-\tilde v_i(k)\big)+a_{0}\big(v_{ix}(k) -\tilde v_{ix}(k)\big)^2 \nonumber
\\
+c_{ik}\big(v_{i}(k) -\tilde v_{i}(k)\big)^2 \Big] \nonumber 
\\
\leq \sum\limits_{i=0}^{m-1}-b_{ik}\big(v_{ix}(k) -\tilde v_{ix}(k)\big)\big(v_{i}(k) -\tilde v_{i}(k)\big) \nonumber
\\
\leq \sum\limits_{i=0}^{m-1} \frac{a_0}{2}\big(v_{ix}(k) -\tilde v_{ix}(k)\big)^2 + \sum\limits_{i=0}^{m-1}  \frac{\|b\|^2_{L^\infty}}{2a_0}\big(v_{i}(k) -\tilde v_{i}(k)\big)^2
\end{gather}

Absorbing to left hand side, and by using \eqref{b_n bound}, we get:

\begin{gather}
\sum\limits_{i=0}^{m-1}\Big[\Big(\frac{\bar{b}}{\tau}-\|c\|_{L^\infty}-\frac{\|b\|^2_{L^\infty}}{2a_0}\Big) \big(v_i(k)-\tilde v_i(k)\big)^2+\frac{a_0}{2}\big(v_{ix}(k) -\tilde v_{ix}(k)\big)^2
 \Big] \leq 0
\end{gather}

The whole summand is non-negative for sufficiently small $\tau$. Therefore, it is equal to 0, which  implies that $v_i(k)=\tilde v_i(k),\quad \forall i=\overline{0,m}$. Hence, by induction, $v=\tilde v$.

Now we seek to prove existence through induction. Construct $v(0)$ through definition of a Discrete State Vector. Note that $v(0)$ is bounded since $\Vert v(0)\Vert \leq  \Vert\Phi\Vert_{L_{\infty}[0,\ell]}$.  Fix $k\geq 1$, and assume that $v(k-1)$ has been constructed so that (\ref{dsvsum}) is satisfied for all $K<k$. Moreover, assume that each element of $v(k-1)$ is bounded. Through manipulation, the summation identity (\ref{dsvsum}) is equivalent to solving the following system of non-linear equations:
\begin{gather}
[a_{0k}-hb_{0k}+h^2c_{0k}]v_0(k)+\frac{h^2}{\tau}b_n(v_0(k))-a_{0k}v_1(k)\nonumber 
\\
=\frac{h^2}{\tau}b_n(v_0(k-1))+h^2f_{0k}-hg_0^n\nonumber
\\ \nonumber \\
\frac{h^2}{\tau}b_n(v_i(k))+[-a_{i-1,k}+hb_{i-1,k}]v_{i-1}(k)+[a_{i-1,k}+a_{ik}-hb_{ik}+h^2c_{ik}]v_i(k) \nonumber
\\ 
-a_{ik}v_{i+1}(k) = \frac{h^2}{\tau}b_n(v_i(k-1))+h^2f_{ik},\quad i =\overline{1,m-1}\nonumber
\\ \nonumber \\
[-a_{m-1,k}+hb_{m-1,k}]v_{m-1}(k)+a_{m-1,k}v_m(k)=hp_k \label{system}
\end{gather}

We will construct $v(k)$ by the method of successive approximations. Fix $h$ and $\tau$, and choose $v^0=v(k-1)$. Having obtained $v^N$, we search $v^{N+1}$ as a solution of the following:
\begin{gather}
[a_{0k}-hb_{0k}+h^2c_{0k}]v_0^{N+1}(k)+\frac{h^2}{\tau}b_n(v_0^{N+1}(k))-a_{0k}v_1^N(k)\nonumber 
\\
=\frac{h^2}{\tau}b_n(v_0(k-1))+h^2f_{0k}-hg_0^n \nonumber
\\ \nonumber \\
\frac{h^2}{\tau}b_n(v^{N+1}_i(k))+[-a_{i-1,k}+hb_{i-1,k}]v_{i-1}^N(k)+[a_{i-1,k}+a_{ik}-hb_{ik}\nonumber
\\
+h^2c_{ik}]v_i^{N+1}(k)-a_{ik}v_{i+1}^N(k) \nonumber 
\\
 = \frac{h^2}{\tau}b_n(v_i(k-1))+h^2f_{ik},\quad i =\overline{1,m-1} \nonumber 
\\ \nonumber \\
 [-a_{m-1,k}+hb_{m-1,k}]v_{m-1}^{N+1}(k)+a_{m-1,k}v_m^{N+1}(k)=hp_k \label{system1}
\end{gather}

We now proceed to prove that the sequence $\{v^N\}$ converges to the unique solution of \eqref{system}.
Subtract \eqref{system1} for $N$ and $N-1$ to get
\begin{gather}
[a_{0k}-hb_{0k}+h^2c_{0k}]\big(v_0^{N+1}(k)-v_0^N(k)\big)+\frac{h^2}{\tau}\big(b_n(v_0^{N+1}(k))-b_n(v_0^{N}(k))\big)
\\
=a_{0k}\big(v_1^N(k)-v_1^{N-1}(k)\big) \nonumber
\\ \nonumber \\
[a_{ik}+a_{i-1,k}-hb_{ik}+h^2c_{ik}]\big(v_i^{N+1}(k)-v_i^N(k)\big)+\frac{h^2}{\tau}\big(b_n(v_i^{N+1}(k))-b_n(v_i^{N}(k))\big)\nonumber
\\
=a_{ik}\big(v_{i+1}^N(k)-v_{i+1}^{N-1}(k)\big) + [a_{i-1,k}-hb_{i-1,k}]\big(v_{i-1}^N(k)-v_{i-1}^{N-1}(k)\big),\quad i =\overline{1,m-1} \nonumber
\\ \nonumber \\ 
a_{m-1,k}\Big(v_m^{N+1}(k)-v_m^N(k)\Big)=[a_{m-1,k}-hb_{m-1,k}]\Big(v_{m-1}^{N+1}(k)-v_{m-1}^{N}(k)\Big)
\end{gather}

which can be transformed to
\begin{gather}\label{system2}
\begin{cases}
v_0^{N+1}(k)- v_0^N(k)&=\Big(\frac{a_{0k}}{a_{0k}-hb_{0k}+h^2c_{0k}+\frac{h^2}{\tau}\zeta^0_{n,N}}\Big)\big(v_1^N(k)-v_1^{N-1}(k)\big)\\[4mm]
v_i^{N+1}(k)-v_i^N(k)&=\Big(\frac{a_{ik}}{a_{i-1,k}+a_{ik}-hb_{ik}+h^2c_{ik}+\frac{h^2}{\tau}\zeta^i_{n,N}}\Big)(v_{i+1}^N(k)- v_{i+1}^{N-1}(k))
\\
&\quad +\Big(\frac{a_{i-1,k}-hb_{i-1,k}}{a_{i-1,k}+a_{ik}-hb_{ik}+h^2c_{ik}+\frac{h^2}{\tau}\zeta^i_{n,N}}\Big)(v_{i-1}^N(k)-v_{i-1}^{N-1}(k))\\[4mm]
v_m^{N+1}(k)-v_m^N(k)&= \frac{a_{m-1,k}-hb_{m-1,k}}{a_{m-1,k}}\Big(v_{m-1}^{N+1}(k)-v_{m-1}^N(k)\Big)
\end{cases}
\end{gather}
where 
\begin{equation*}
\zeta_{n,N}^i:=\int_0^1b_n'(\theta v_i^{N+1}(k)+(1-\theta)v_i^N(k))d\theta, \quad i =\overline{0,m-1}.
\end{equation*}
Due to \eqref{b_n bound}, we have $\zeta_{n,N}^i\geq\bar b,~~i=\overline{0,m-1}$.
Let 
\[A_N:=\max\limits_{0\leq i\leq m} {\big|v_i^{N+1}(k)-v_i^N(k)\big|}. \] 
From (\ref{system2}), taking the first equation into consideration, we have:
\begin{gather}\label{delta0}
|v_0^{N+1}(k)- v_0^N(k)| \leq \Big | \frac{a_{0k}}{a_{0k}-hb_{0k}+h^2c_{0k}+\frac{h^2}{\tau}\zeta^0_{n,N}}\Big | A_{N-1} \nonumber \\
=|\delta_0^{-1}|A_{N-1}, \ \delta_0=1+\frac{h(-b_{0k}+hc_{0k}+\frac{h}{\tau}\zeta^0_{n,N})}{a_{0k}}
\end{gather}
We have
\[ 0<a_0\leq a_{0k} \leq \|a\|_{L_\infty} \]
and
\begin{gather*}
	-b_{0k}+hc_{0k}+\frac{h}{\tau}\zeta^0_{n,N} \geq -\|b\|_{L_\infty}-h\|c\|_{L_\infty} +\frac{h}{\tau}\bar{b} > 0 
\end{gather*}
by \eqref{htau} and for sufficiently small $h$ and $\tau$. Thus, $0<\delta_0^{-1} <1$. Similarly,
\begin{gather}\label{deltai}
|v_i^{N+1}(k)- v_i^N(k)| \leq  \Big|\frac{a_{i-1,k}+a_{ik}-hb_{ik}}{a_{i-1,k}+a_{ik}-hb_{ik}+h^2c_{ik}+\frac{h^2}{\tau}\zeta^i_{n,N}}\Big|A_{N-1}\nonumber\\
=|\delta_i^{-1}|A_{N-1}, \ \delta_i=1+\frac{h(hc_{ik}+\frac{h}{\tau}\zeta^i_{n,N})}{a_{i-1,k}+a_{ik}-hb_{ik}}, i = \overline{1,m-1}
\end{gather}
Through similar argument as with $\delta_0$, we derive that $0<\delta_i^{-1} <1$ for $i = \overline{1,m-1}$ for $0<h<<1$. For $i=m$, we get

\begin{gather*}\label{deltam}
|v_m^{N+1}(k)- v_m^N(k)|\leq |\delta_m^{-1}| A_{N-1}, \\
\delta_m= \Big(1-\frac{hb_{m-1,k}}{a_{m-1,k}}\Big)^{-1}\Big(1+\frac{h^2c_{m-1,k}+\frac{h^2}{\tau}\zeta^{m-1}_{n,N}}{a_{m-2,k}+a_{m-1,k}-hb_{m-1,k}}\Big)
\end{gather*}

For $0<h<<1$, we can see that the term in left brackets will be close to 1, and due to \eqref{htau}, as before we derive that $0<\delta_m^{-1}<1$ for sufficiently small $h$ and $\tau$. Let $\delta=\max\limits_{i=\overline{1,m}}\delta_i^{-1}$. We thus have $\delta<1$, and 

\begin{equation}\label{A}
A_N\leq\delta A_{N-1}\leq\cdots\leq A_0\delta^N.
\end{equation}
Following the proof given in \cite{Abdulla5} (Lemma 1, Section 2) it follows that there exist finite limits
\begin{equation}\label{limit}
v_i(k)=\lim\limits_{N\rightarrow+\infty}v_i^N(k),\qquad i=0,1,\ldots,m.
\end{equation}
Passing to limit as $N\to+\infty$ in \eqref{system1}, we derive that $v_i(k), i=\overline{1,m}$ is a unique solution of \eqref{system}.\hfill{$\square$}

Given the existence and uniqueness of the discrete state vector for fixed $n$, we can uniquely define for each $k=1,\ldots,n$ the vector $\zeta_k$ whose $m$ components $\zeta_k^i$ are given by
\begin{equation}
\zeta^i_k=\int_0^1b_n'(\theta v_i(k)+(1-\theta)v_i(k-1))d\theta, \quad i =\overline{0,m-1}\label{zetav}.
\end{equation}

The following is a well known necessary and sufficient condition for the convergence of the discrete optimal control problems to continuous optimal control problem.

\begin{lemma}\label{Vasil}\cite{Vasilev} The sequence of discrete optimal control problems approximates the continuous optimal control problem if and only if the following conditions are satisfied:
	\begin{itemize}
		\item for arbitrary sufficiently small $\ep>0$ there exists $M_1=M_1(\ep)$ such that $\n{Q}_M(g)\in \n G^{M}_R$ for all $g \in \n G_{R-\ep}$ and $M\ge M_1$; and for any fixed $\ep>0$ and for all $g\in \n G_{R-\ep}$ the following inequality is satisfied:
		\begin{equation}\label{firstcondition}
		\limsup\limits_{M\to \infty} \Big ( \n{I}_M(\n{Q}_M(g))-\n{J}(g) \Big ) \le 0.
		\end{equation}
		\item for arbitrary sufficiently small $\ep>0$ there exists  $M_2=M_2(\ep)$ such that $\n{P}_M([g]_{M})\in \n G_{R+\ep}$ for all $[g]_{M} \in \n G^{M}_R$ and $M\ge M_2$; and for all $[g]_{M}\in \n G^{M}_R$, $M\ge 1$ the following inequality is satisfied:
		\begin{equation}\label{secondcondition}
		\limsup\limits_{M\to \infty} \Big ( \n{J}(\n{P}_M([g]_{M})) -\n{I}_M([g]_{M})  \Big ) \le 0.
		\end{equation}
		\item the following inequalities are satisfied:
		\begin{align}
		\limsup\limits_{\ep \to 0} \n{J}_*(\ep) \ge \n{J}_*, 
		\ \ \liminf\limits_{\ep \to 0} \n{J}_*(-\ep) \le \n{J}_*,
		\end{align}\label{thirdcondition}
		where $\n{J}_*(\pm\ep)=\inf\limits_{\n G_{R\pm \ep}}\n{J}(g)$.
	\end{itemize}
\end{lemma}

\begin{lemma}\label{mappings}\cite{Abdulla5} The mappings $\n P_n, \n Q_n$ satisfy the conditions of Lemma \ref{Vasil}.
\end{lemma}

\begin{lemma}\label{unique} There is at most one solution to the multiphase free boundary problem (\ref{bvpde})-(\ref{vp}) in the sense of (\ref{weaksolution}). 
\end{lemma}

\noindent\emph{Proof.} The uniqueness of the weak solution is proved in Section $9$ of Chapter V of \cite{LSU} for the classical multiphase Stefan Problem 
($\mathcal{L}=\Delta$) under zero Dirichlet boundary conditions on the fixed boundary. Lemma 4 of \cite{Abdulla5} generalized the result to the case of non-homogeneous Neumann boundary condition. We generalize the result to the case of multiphase free boundary problem with general elliptic operator $\mathcal{L}$. Uniqueness is proved over a wider class of solutions than given in (\ref{weaksolution}). Suppose that $v\in L_{\infty}(D)$ only, not necessarily in the Sobolev space $W_2^{1,1}(D)$, and that for any two functions $B,B_0$ of type $\n B$ it satisfies the identity
\begin{gather}
\int\limits_0^T\int\limits_0^{\ell}\Big[B(x,t,v)\psi_t+v(a\psi_{x})_x-bv \psi_x-cv\psi+f\psi\Big]\,dxdt \nonumber \\ + \int\limits_0^{\ell}B_0(x,0,\Phi(x))\psi(x,0)\,dx  +\int\limits_0^Tp(t)\psi(\ell,t)\,dt-\int\limits_0^Tg(t)\psi(0,t)\,dt=0, \nonumber 
\\
\forall \psi\in W_2^{2,1}(D), \psi(x,T)=0, \nonumber
\\
 a(0,t)\psi_x(0,t)=a(\ell,t)\psi_x(\ell,t)=0. \label{weakersol}
\end{gather}
Any function satisfying (\ref{weaksol}) will also satisfy the above definition. Suppose $v$ and $\tilde v$ are two solutions in the sense of (\ref{weakersol}), and subtract (\ref{weakersol}) with solution $\tilde{v}$ from that of $v$. Due to $\Phi$ taking on phase transition temperatures on sets of measure 0, the $B_0$ term will vanish and we are left with the following:
\begin{gather}\label{psiintegral}
\int\limits_0^T\int\limits_0^{\ell}\big(B(x,t,v)-\tilde B(x,t,\tilde v)\big)\left(\psi_t+z(x,t)\Big((a\psi_{x})_x-b\psi_x-c\psi\Big) \right)\,dx\,dt = 0
\end{gather}
where $z(x,t) = \frac{v-\tilde v}{B(x,t,v)-\tilde B(x,t,\tilde v)}$. For $(x,t)\in D$ such that $v(x,t)=\tilde v(x,t)$, we have $z(x,t)=0$. Otherwise, since $B$ and $\tilde B$ are strictly increasing on $v$ a.e. $(x,t)\in D$, we have that $z$ is non-negative for a.e. $(x,t)$. Moreover, we have:
\begin{align*}
|z(x,t)| & = \left|\frac{v-\tilde v}{\int\limits_{\tilde v(x,t)}^{v(x,t)}\beta'(w)\,dw+\sum\limits_{i: v^i\in(\tilde{v}(x,t), v(x,t))}(\beta(v^i)^+-\beta(v^i)^-)}\right| \leq\left|\frac{v-\tilde v}{\int\limits_{\tilde v}^v\bar b\,dv}\right| = \frac1{\bar b},
\end{align*}
so that $z$ is essentially bounded. Fix $\ep>0$, and take as $\psi(x,t)$ the solution of the following Neumann problem
\begin{gather}
\psi_t+(z(x,t)+\ep)\Big((a\psi_{x})_x-b\psi_x-c\psi\Big)=F(x,t)\label{conjugate},
\\
a(0,t)\psi_x(0,t)=a(\ell,t)\psi_x(\ell,t)=0. \label{psiboundary}
 \\
 \quad \psi(x,T)=0\label{psifinal},
\end{gather}
where the $\ep$ is added to ensure the conjugate diffusion coefficient is strictly positive, and $F$ is an arbitrary smooth bounded function in $D$. Note that (\ref{conjugate}) is the conjugate parabolic equation. From \cite{LSU}, there exists a unique solution $\psi^{\ep} \in W_2^{2,1}(D)$ of the problem \eqref{conjugate}-\eqref{psifinal}. We will use the arbitrariness of $F$ to obtain that $B-\tilde B=0$ a.e. $(x,t) \in D$. Note that through \eqref{conjugate}, we can rewrite \eqref{psiintegral}:
\begin{equation}
\int\limits_0^T\int\limits_0^{\ell}\big(B(x,t,v)-\tilde B(x,t,\tilde v)\big)\left(F-\ep\Big((a\psi_{x})_x-b\psi_x-c\psi\Big) \right)\,dx\,dt = 0 \label{almost}.
\end{equation}
Thus our goal will be attained if we have an energy estimate on $\mathcal{L}\psi$ for solutions of (\ref{conjugate}). For simplicity, we obtain energy estimates through the second order parabolic equation, which will give analogous estimates for the conjugate parabolic equation by reversing the time variable. Let $z^{\ep}(x,t) =z(x,t)+\ep$, and for simplicity we will omit the superscript. Multiply the non-conjugate version of (\ref{conjugate}) by $\psi_{xx}$ and integrate it over $D_t :=(0,\ell)\times(0,t)$ to get
\begin{gather}
-\int\limits_0^t\int\limits_0^{\ell}(\psi_{\tau}-za\psi_{xx}-z(a_x-b)\psi_x+zc\psi)\psi_{xx}\,dx\,d\tau = -\int\limits_0^t\int\limits_0^{\ell}F\psi_{xx} \,dx\,d\tau \nonumber\\
= \int\limits_0^t\int\limits_0^{\ell}F_x\psi_x\,dx\,d\tau-\int\limits_0^tF\psi_x\Big|_0^{\ell}\,d\tau, \label{psixxtest} 
\end{gather}
Due to \eqref{psiboundary}, the second integral on the right hand side disappears. We can transform various terms on the right hand sign as follows:
\begin{gather*}
-\int\limits_0^t\int\limits_0^{\ell} \psi_{\tau}\psi_{xx}\,dx\,d\tau = \int\limits_0^t\int\limits_0^{\ell} (\psi_{\tau})_x\psi_{x}\,dx\,d\tau-\int\limits_0^t \psi_{\tau}\psi_{x}\Big|_0^{\ell} \,d\tau = \frac{1}{2}\int\limits_0^{\ell}\psi_x^2(x,\ell)\,dx, \\
- \int\limits_0^t\int\limits_0^{\ell} (-za\psi_{xx})\psi_{xx}\,dx\,d\tau \geq a_0\int\limits_0^t\int\limits_0^{\ell} z\psi_{xx}^2 \,dx\,d\tau
\end{gather*}
Using the above, and returning to \eqref{psixxtest}, we get that:
\begin{gather}
\frac{1}{2}\int_0^{\ell}\psi_x^2(x,\ell)\,dx+a_0\int\limits_0^t\int\limits_0^{\ell} z\psi_{xx}^2 \,dx\,d\tau \leq \int\limits_0^t\int\limits_0^{\ell} (-za_x\psi_x\psi_{xx}) \,dx\,d\tau \nonumber\\
+ \int\limits_0^t\int\limits_0^{\ell} zb\psi_x\psi_{xx} \,dx\,d\tau+\int\limits_0^t\int\limits_0^{\ell} zc\psi\psi_{xx} \,dx\,d\tau + \int\limits_0^t\int\limits_0^{\ell}F_x\psi_x\,dx\,d\tau \label{psixxtest2}
\end{gather}
We now estimate the terms on the right hand side using Cauchy inequality with $\epsilon>0$ and properties of given functions,
and absorbing terms to the left hand side, we have:
\begin{gather}
\frac{1}{2}\int_0^{\ell}\psi_x^2(x,\ell)\,dx+\frac{a_0}{4}\int\limits_0^t\int\limits_0^{\ell} z\psi_{xx}^2 \,dx\,d\tau \nonumber \\
\leq \Big(\frac{2(\Vert a_x \Vert^2_{L_{\infty}(D)}+\Vert b \Vert^2_{L_{\infty}(D)})}{\bar{b}a_0}+\frac{1}{2}\Big)\int\limits_0^t\int\limits_0^{\ell} \psi_x^2 \,dx\,d\tau +\frac{1}{2} \int\limits_0^t\int\limits_0^{\ell}F^2_x\,dx\,d\tau\nonumber \\
+ \frac{2\Vert c \Vert^2_{L_{\infty}(D)}}{\bar{b}a_0}\int\limits_0^t\int\limits_0^{\ell}\psi^2 \,dx\,d\tau, \label{psixxtest3}	
\end{gather}
From Theorem 2.3. Chapter 1 of \cite{LSU} it follows that $\psi$ will have a uniform bound in $L_\infty(D)$, which we will denote as $\bar{C}$. Letting now $y(t) = \int\limits_0^t\int\limits_0^{\ell}\psi_x^2\,dx\,d\tau$ from (\ref{psixxtest3}) we deduce
\[
y'(t)\leq Cy(t)+\Big(\int\limits_0^t\int\limits_0^{\ell}F_x^2\,dx\,d\tau+\frac{4\Vert c \Vert^2_{L_{\infty}(D)}\bar{C}^2\ell T}{\bar{b}a_0}\Big).
\]
where 
\[ C = 2\Big(\frac{2(\Vert a_x \Vert^2_{L_{\infty}(D)}+\Vert b \Vert^2_{L_{\infty}(D)})}{\bar{b}a_0}+\frac{1}{2}\Big)\] 
By Gronwall's Inequality (e.g. Lemma 5.5, Chapter 2, \cite{LSU}), we deduce from the above differential inequality that 
\begin{gather*}
\int\limits_0^t\int\limits_0^{\ell}\psi_x^2(x,\tau) \,dx\,d\tau \leq \Big[\frac{e^{Ct}-1}{C}\Big]\Big[\int\limits_0^t\int\limits_0^{\ell}F_x^2\,dx\,d\tau+\frac{4\Vert c \Vert^2_{L_{\infty}(D)}\bar{C}^2\ell T}{\bar{b}a_0}\Big], \quad \forall t\in(0,T]
\end{gather*}
so that by (\ref{psixxtest3}),
\begin{gather*}
\int\limits_0^{\ell}\psi_x^2(x,t) \,dx + \frac{a_0}{2}\int\limits_0^t\int\limits_0^{\ell}z\psi_{xx}^2 \,dx\,d\tau \leq (e^{Ct}-1))\Big[\int\limits_0^t\int\limits_0^{\ell}F_x^2\,dx\,d\tau\\
+\frac{4\Vert c \Vert^2_{L_{\infty}(D)}\bar{C}^2\ell T}{\bar{b}a_0}\Big]+\frac{1}{2}\int\limits_0^t\int\limits_0^{\ell}F_x^2\,dx\,d\tau+\frac{4\Vert c \Vert^2_{L_{\infty}(D)}\bar{C}^2\ell T}{\bar{b}a_0}\, \qquad \forall t\in(0,T],
\end{gather*}
The first of the above inequalities implies that 
\begin{gather*}
\underset{0\leq t\leq T}{\text{ess sup}}\int\limits_0^{\ell}\psi_x^2(x,t) \,dx\leq C_1\int\limits_0^T\int\limits_0^{\ell}F_x^2\,dx\,d\tau+C_2. 
\end{gather*}
Now, since $\psi_t = az\psi_{xx}+z(a_x-b)\psi_x+cz\psi+F$, we have
\begin{gather*}
\Vert\psi_t\Vert^2_{L_2(D_t)} = \Vert az\psi_{xx}+z(a_x-b)\psi_x+cz\psi +F\Vert^2_{L_2(D_t)} \\
\leq 4\Big(\Vert az\psi_{xx}\Vert^2_{L_2(D_t)}+\Vert (a_x-b)z\psi_{x}\Vert^2_{L_2(D_t)}+\Vert cz\psi\Vert^2_{L_2(D_t)}+\Vert F\Vert^2_{L_2(D_t)}\Big )\\
 \leq C_3\Big[\Vert F_x\Vert^2_{L_2(D_t)}+\Vert F\Vert^2_{L_2(D_t)}\Big] + C_4,
\end{gather*}
where the constants $C_3$ and $C_4$ depend on $\bar{b}, a_0, T, \ell$ and $L_\infty$-norms of $a, a_x, b, c$. Combining all the estimations we have the following desired energy estimate for $\psi\in W_2^{2,1}(D)$:
\begin{gather}
\int\limits_0^T\int\limits_0^\ell(\psi_t^2+z(x,t)\psi_{xx}^2)\,dx\,dt+\underset{0\leq t\leq T}{\text{ess sup}}\int\limits_0^{\ell}\psi_x^2(x,t) \,dx\nonumber\\
\leq C_5\Big[\Vert F_x\Vert^2_{L_2(D_t)}+\Vert F\Vert^2_{L_2(D_t)}\Big] + C_6,\label{energy}
\end{gather}
where the constants $C_5$ and $C_6$ are independent of $\epsilon$, and  depend on $\bar{b}, a_0, T, \ell$ and $L_\infty$-norms of $a, a_x, b, c$.

For the rest of the proof, any constant depending on the bounded data, domain, $F$, or the uniform bound on $\psi$ will be referred to as $C^*$.We can now observe that
\begin{gather*}
\left|\int\limits_0^T\int\limits_0^{\ell}(B-\tilde B)\ep\Big((a\psi^{\ep}_x)_{x}-b\psi^{\ep}_x-c\psi^{\ep}\Big)\,dx\,dt\right|  \\
\leq  \int\limits_0^T\int\limits_0^{\ell}\left|B-\tilde B\right|\ep\Big(\left|a\psi^{\ep}_{xx}\right|+\left|(a_x-b)\psi^{\ep}_x\right|+\left|c\psi^\ep\right|\Big)\,dx\,dt\\
 \leq 2\text{esssup}~\beta(v)\Bigg[ \int\limits_0^T\int\limits_0^{\ell} \frac{\ep(z+\ep)^{\frac12}}{(z+\ep)^{\frac12}}\left|a\psi^{\ep}_{xx}\right|\,dx\,dt \\
 + \int\limits_0^T\int\limits_0^{\ell} \ep \left|a_x-b\right|\left|\psi_x^{\ep}\right|\,dx\,dt + \int\limits_0^T\int\limits_0^{\ell} \ep \left|c\right|\left|\psi^{\ep}\right|\,dx\,dt\Bigg] \\
 \leq 2\text{esssup}~\beta(v)\Bigg[ \Big(\int\limits_0^T\int\limits_0^{\ell} \frac{\ep^2}{z+\ep} \,dx\,dt \Big)^{\frac12}\Big(a^2(z+\ep)(\psi_{xx}^{\ep})^2 \,dx\,dt\Big)^{\frac12} \\
 + \Big(\int\limits_0^T\int\limits_0^{\ell} \ep^2(a_x-b)^2\,dx\,dt \Big)^{\frac12}\Big( \int\limits_0^T\int\limits_0^{\ell} (\psi_x^\ep)^2 \,dx\,dt\Big)^{\frac12}+\ep C^*\Bigg] \\
 \leq 2\sqrt{\ep}~\text{esssup}~\beta(v) \Bigg[ C^*\Big(\int\limits_0^T\int\limits_0^{\ell} \frac{\ep}{z+\ep} \,dx\,dt \Big)^{\frac12}  + C^*\sqrt{\ep} \Bigg]\\
 \leq 2C^*\sqrt{\ep}~\text{esssup}~\beta(v)[(T\ell)^{\frac{1}{2}}+\sqrt{\ep}]\to 0,
\end{gather*} 
 as $\ep \to 0$. Therefore, (\ref{almost}) now implies
\[
\int\limits_0^T\int\limits_0^{\ell}\big(B(x,t,v(x,t))-\tilde B(x,t,\tilde v(x,t))\big)F\,dx\,dt = 0.
\]
Since $F$ is arbitrary, the above equality gives that $B(x,t,v(x,t)) = \tilde B(x,t,\tilde v(x,t))$ a.e. $(x,t)\in D$. This implies $\beta(v(x,t))=\beta(\tilde v(x,t))$, a.e.\,$(x,t) \in D~$s.t.$~v(x,t)\neq v^j, j=1,...,m$. Due to the fact that $\beta$ is strictly increasing, we have $v(x,t)=\tilde v(x,t)$ a.e.\,$(x,t)$. Thus $v$ and $\tilde v$ are the same solution to \eqref{weakersol} \hfill{$\square$}

\begin{corollary}\label{Cjmeasure0} If $v$ is weak solution, the sets $\{(x,t) \in D \vert v=v^j\}, j=1,...,m$ have 2-dimensional measure 0.
\end{corollary}
Indeed, due to uniqueness of the weak solution, for any two representatives $B_1, B_2$ of the class $\n B$ we have
\[ B_1(x,t,v(x,t))=B_2(x,t,v(x,t)), \ \quad\text{a.e.} \ (x,t)\in D \]
and by the Definition~\ref{typeB} this will be a contradiction if any of the $v^j$-level sets of the weak solution would have a positive 2-dimensional measure.
\section{Proof of Main Results}

\subsection{$L_\infty(D)$ estimate for the discrete multiphase free boundary problem}

In this section we prove $L_\infty(D)$ bound for the discrete PDE problem under the following reduced assumptions:
\[  p\in L_{\infty}(0,T),\  \Phi\in L_{\infty}(0,\ell),\  f,a,b,c\in L_{\infty}(D) \]
and $a$ satisfies \eqref{apositive}.
\begin{theorem}\label{uniformboundedness} For $[g]_n\in \n G_R^n$ and $n,m$ large enough, the discrete state vector $[v([g]_n)]_n$ satisfies the following estimate:
	\begin{gather}\label{linfestimate}
	\Vert[v]_n\Vert_{\ell_{\infty}}:=\max\limits_{0\leq k\leq n}\Big(\max\limits_{0\leq i\leq m}|v_i(k)|\Big)\nonumber \\ \leq C_{\infty}\Big(\Vert f\Vert_{L_{\infty}(D)}+\Vert p\Vert_{L_{\infty}(0,T)}+\Vert g^n\Vert_{L_\infty(0,T)}+\Vert\Phi\Vert_{L_{\infty}(0,\ell)}\Big)
	\end{gather}
	where $C_{\infty}$ is a constant independent of $n$ and $m$.
\end{theorem}
\noindent\emph{Proof. } Fix $n$ arbitrarily large. Note $\max|v_i(0)|\leq\Vert\Phi\Vert_{L_{\infty}(0,\ell)}$. Consider a positive function $\gamma(x)\in C^2[0,\ell]$ satisfying
\begin{gather}\label{gamma}
\gamma(0)=\frac12,\ \gamma(\ell)=\frac12,\ \gamma'(0)=\frac{4(\Vert b\Vert_{L_\infty(D)}+\Vert c\Vert_{L_\infty(D)})}{a_0}+1,\nonumber \\ \gamma'(\ell)=-1,\qquad \frac14\leq\gamma(x)\leq1,~x\in[0,\ell].
\end{gather}
Define $\gamma_i=\gamma(x_i),~i=\overline{0,m}$, and denote as $x^i$ the value in $[x_i,x_{i+1}]$ that satisfies (by mean value theorem (MVT)) $\gamma(x_{i+1})-\gamma(x_i)=\gamma'(x^i)h$. Transform the discrete state vector as
\[
w_i(k)=v_i(k)\gamma_i,\qquad i =\overline{0,m},\qquad k =\overline{0,n}.
\]
System (\ref{system}) can be rewritten as:
\begin{gather}\label{dpde}
\left\{\begin{matrix}
	h\zeta_k^0v_{0\bar t}(k)-a_{0k}v_{0x}(k)+[hc_{0k}-b_{0k}]v_{0}(k)= hf_{0k}-g_k^n&\\[4mm]
\zeta_k^iv_{i\bar t}(k)-[a_{i-1,k}-hb_{ik}]v_{ix\bar x}(k)+\frac{1}{h}[a_{i-1,k}-hb_{ik}-a_{ik}]v_{ix}(k)\\+c_{ik}v_i(k) = f_{ik}&,\, i =\overline{1,m-1} \\[4mm] 
a_{m-1,k}v_{m-1,x}(k)+b_{m-1,k}v_{m-1}(k)=p_k& \end{matrix}\right..
\end{gather}
We note
\begin{gather*}
v_i(k) = \frac1{\gamma_i}w_i(k),\qquad v_{i\bar t}(k) = \frac1{\gamma_i}w_{i\bar t}(k), \\[4mm]
v_{ix}(k) = \frac1{\gamma_{i+1}}w_{ix}(k)+\left(\frac1{\gamma_i}\right)_xw_i(k) = \frac1{\gamma_{i}}w_{ix}(k)+\left(\frac1{\gamma_i}\right)_xw_{i+1}(k), \\[4mm]
v_{ix\bar x}(k) = \frac1{\gamma_{i-1}}w_{ix\bar x} (k) +\left[ \left(\frac1{\gamma_i}\right)_{\bar x} + \left(\frac1{\gamma_i}\right)_{x} \right]w_{ix}(k)+ \left(\frac1{\gamma_i}\right)_{x\bar x}w_i(k) \\[2mm]\qquad\qquad\qquad =\frac1{\gamma_{i+1}}w_{ix\bar x} (k) +\left[ \left(\frac1{\gamma_i}\right)_{\bar x} + \left(\frac1{\gamma_i}\right)_{x} \right]w_{i\bar x}(k)+ \left(\frac1{\gamma_i}\right)_{x\bar x}w_i(k), \\[4mm]
\left(\frac1{\gamma_i}\right)_x = -\frac1{\gamma_i\gamma_{i+1}}\gamma_{ix},\qquad \left(\frac1{\gamma_i}\right)_{x\bar x} = -\frac1{\gamma_i\gamma_{i+1}}\gamma_{ix\bar x}+\frac{\gamma_{ix}+\gamma_{i\bar x}}{\gamma_{i-1}\gamma_i\gamma_{i+1}}\gamma_{i\bar x}.
\end{gather*}
Thus $w_i(0)=\gamma_i\Phi_i,\quad i=\overline{0,m}$, and for $k=\overline{1,n}$,
\begin{equation}\label{systemw}
\left\{\begin{matrix}\frac h{\gamma_0}\zeta_k^0w_{0\bar t}(k)-\frac{a_{0k}}{\gamma_1}w_{0x}(k)-\Big[a_{0k}\left(\frac1{\gamma_0}\right)_x+\frac{b_{0k}-hc_{0k}}{\gamma_0}\Big]w_0(k)= hf_{0k}-g_k^n\\ \\[4mm]
\frac1{\gamma_i}\zeta_k^iw_{i\bar t}(k)-\frac{a_{i-1,k}-hb_{ik}}{\gamma_{i-1}}w_{ix\bar x}(k) \\-\left[ (a_{i-1,k}-hb_{ik})\Big(\left(\frac1{\gamma_i}\right)_{\bar x} + \left(\frac1{\gamma_i}\right)_{x}\Big)+\frac1{\gamma_{i+1}h}[-a_{i-1,k}+hb_{ik}+a_{ik}] \right]w_{ix}(k) \\-\left[ (a_{i-1,k}-hb_{ik})\left(\frac1{\gamma_i}\right)_{x\bar x} +\frac1{h}(-a_{i-1,k}+hb_{ik}+a_{ik})\left(\frac1{\gamma_i}\right)_x+\frac{c_{ik}}{\gamma_i}\right]w_i(k)\\=f_{ik},\quad i =\overline{1,m-1} \\ \\[4mm] 
\frac{a_{m-1,k}-b_{m-1,k}h}{\gamma_{m-1}}w_{m-1,x}(k)+\left[a_{m-1,k}\left(\frac{1}{\gamma_{m-1}}\right)_x+\frac{b_{m-1,k}h}{\gamma_{m-1}}\right]w_m(k)=p_k \end{matrix}. \right..
\end{equation}
Furthermore, transform $w_i(k)$ as:
\begin{equation}
\label{u} u_i(k) = w_i(k)e^{-\lambda t_k},\qquad i = \overline{0,m},\quad k =\overline{0,n}
\end{equation}
where $\lambda$ satisfies 
\begin{gather}\label{lambda}
\bar{b}(\lambda-1) = (\|a\|_{L_\infty(D)}+\|b\|_{L_\infty(D)})(32\|\gamma''\|_{C[0,\ell]}+356\|\gamma'\|^2_{C[0,\ell]})\nonumber \\
+32(\|a_x\|_{L_\infty(D)}+\|b\|_{L_\infty(D)})\|\gamma'\|_{C[0,\ell]}+8\|c\|_{L_\infty(D)}
\end{gather}
and if $t^k\in[t_{k-1},t_k]$ satisfies through the MVT that $e^{\lambda t_k}-e^{\lambda t_{k-1}} = \lambda e^{\lambda t^k}\tau$, then 
\[
w_{i\bar t}(k) = e^{\lambda t_{k-1}}u_{i\bar t}(k)+\lambda e^{\lambda t^k}u_i(k).
\] 
So $u_i(0) = w_i(0)=\gamma_i\Phi_i,\quad i=\overline{0,m}$, and for $k=\overline{1,n}$, the vector $u(k)$ satisfies the system
\begin{gather}
\frac{h}{\gamma_0}\zeta_k^0e^{-\lambda \tau}u_{0\bar t}(k)-\frac{a_{0k}}{\gamma_1}u_{0x}(k)+\Big[-a_{0k}\left(\frac1{\gamma_0}\right)_x-\frac{b_{0k}-hc_{0k}}{\gamma_0}\nonumber\\+\frac{h\lambda}{\gamma_0}\zeta_k^0e^{-\lambda(t_k-t^k)}\Big]u_0(k)\nonumber = e^{-\lambda t_k}(hf_{0k}-g_k^n), \nonumber\\ \nonumber\\[4mm]
\frac1{\gamma_i}\zeta_k^i e^{-\lambda \tau}u_{i\bar t}(k)-\frac{a_{i-1,k}-hb_{ik}}{\gamma_{i-1}}u_{ix\bar x}(k)\nonumber\\-\left[ (a_{i-1,k}-hb_{ik})\Big(\left(\frac1{\gamma_i}\right)_{\bar x} + \left(\frac1{\gamma_i}\right)_{x}\Big)+\frac1{\gamma_{i+1}h} [-a_{i-1,k}+hb_{ik}+a_{ik}] \right]u_{ix}(k)\nonumber\\-\Bigg[(a_{i-1,k}-hb_{ik})\left(\frac1{\gamma_i}\right)_{x\bar x} +\frac1{h}(-a_{i-1,k}+hb_{ik}+a_{ik})\left(\frac1{\gamma_i}\right)_x \label{systemu}\\+\frac{c_{ik}}{\gamma_i}-\frac{\zeta_k^i\lambda e^{-\lambda(t_k-t^k)}}{\gamma_i}\Bigg]u_i(k)=e^{-\lambda t_k}f_{ik},\quad i =\overline{1,m-1}\nonumber\\  \nonumber\\[4mm]
\frac{a_{m-1,k}-b_{m-1,k}h}{\gamma_{m-1}}u_{m-1,x}(k)+\left[a_{m-1,k}\left(\frac{1}{\gamma_{m-1}}\right)_x+\frac{b_{m-1,k}h}{\gamma_{m-1}}\right]u_m(k)=e^{-\lambda t_k}p_k\nonumber
\end{gather} 
Now fix $k_1\leq n$, and define the following sets of indexes for convenience:
\begin{gather*}
\n M_{k_1} =\{(i,k)|i=0,\ldots,m,\quad k = 0,\ldots, k_1\} , \\
\n N =\{(i,k)|i=1,\ldots,m-1,\quad k = 1,\ldots, k_1\}, \\
\n T_0 = \{(i,k)|i=0,k = 1,\ldots, k_1\}, \\
\n T_m = \{(i,k)|i=m,k = 1,\ldots, k_1\}, \\
\n X_0 = \{(i,k)|i=0,\ldots,m,\quad k =0\}.
\end{gather*}
Unless confusion may arise, we omit the subscript to $\n M_{k_1}$. It is clear that 
\[
\n M = \n N\cup\n T_0\cup\n T_m\cup\n X_0.
\]
If $u_i(k)\leq 0$ in $\n M$, then $\max\limits_{\n M}u_i(k)\leq 0$. Suppose that there exists $(i,k)$ such that $u_i(k)>0$. Then $\max\limits_{\n M}u_i(k)>0$. Let $(i^*,k^*)\in \n M$ be such that $u_{i^*}(k^*) = \max\limits_{\n M}u_i(k)$.

If $(i^*,k^*)\in\n X_0$, then $u_{i^*}(k^*) = \max\limits_i\gamma_i\Phi_i \leq \max\limits_i\Phi_i \leq \max\limits_{[0,\ell]}\Phi(x)$.

If $(i^*,k^*)\in\n T_m$, then $i^*=m,~u_{m-1,x}(k^*)\geq0$ and we can choose $h$ small enough that $\gamma_{m-1,x}=\gamma'(x^{m-1})\in(-\frac32,-\frac12)$ so that
\[
\Bigg(-\frac{a_{m-1,k^*}\gamma'(x^{m-1})}{\gamma_m\gamma_{m-1}}+\frac{b_{m-1,k^*}h}{\gamma_{m-1}}\Bigg)u_m(k^*)\leq e^{-\lambda t_{k^*}}p_{k^*} 
\]

We can see that:

\[
	-\frac{a_{m-1,k^*}\gamma'(x^{m-1})}{\gamma_m\gamma_{m-1}}+\frac{b_{m-1,k^*}h}{\gamma_{m-1}}	 \geq \frac{a_0}{2}-4\Vert b\Vert_{L_\infty(D)}h \geq \frac{a_0}{4}
\]

for sufficiently small h. Thus we have:

\[
	u_m(k^*) \leq \frac{4}{a_0}e^{-\lambda t_{k^*}}p_{k^*} 
\]

If $(i^*,k^*)\in\n T_0$, then $i^*=0, u_{0\bar t}(k^*)\geq0$, $u_{0x}(k^*)\leq0$. Notice that $\left(\frac1{\gamma_0}\right)_x = -\frac1{\gamma_0\gamma_{1}}\gamma_{0x}$. Note $\gamma_{0x}=\gamma'(x^0)$, so for $h$ small enough, we can ascertain $\gamma_{0x}=\gamma'(x^0)\in(\frac{4(\Vert b\Vert_{L_\infty(D)}+\Vert c\Vert_{L_\infty(D)})}{a_0}+\frac12,\frac{4(\Vert b\Vert_{L_\infty(D)}+\Vert c\Vert_{L_\infty(D)})}{a_0}+\frac32)$. It follows
\[
	\Big[a_{0k^*}\frac{\gamma'(x^0)}{\gamma_0\gamma_1}-\frac{b_{0k^*}-hc_{0k^*}}{\gamma_0}+\frac{h\lambda}{\gamma_0}\zeta_{k^*}^0e^{-\lambda(t_{k^*}-t^{k^*})}\Big]u_0(k^*)= e^{-\lambda t_{k^*}}(hf_{0k^*}-g_{k^*}^n)
\]

Since the third term in the parenthesis on the left hand side is positive, we only consider first two terms. We can see that:

\[
	a_{0k^*}\frac{\gamma'(x^0)}{\gamma_0\gamma_1}-\frac{b_{0k^*}-hc_{0k^*}}{\gamma_0} \geq a_0\gamma'(x^0) - 4(\Vert b \Vert_{L_\infty(D)}+h\Vert c \Vert_{L_\infty(D)}) \geq \frac{a_0}{2}
\]

Thus we have:

\[
u_0(k^*) \leq \frac{2}{a_0}	e^{-\lambda t_{k^*}}(hf_{0k^*}-g_{k^*}^n)
\]

If $(i^*,k^*)\in\n N$, then $u_{i^*\bar t}(k^*)\geq 0$, $~u_{i^*x\bar x}(k^*)=\frac1{h^2} \big(u_{i^*+1}(k^*)-2u_{i^*}(k^*)+u_{i^*-1}(k^*)\big)\leq 0$. For $(i,k)\in\n N$, the corresponding equation in (\ref{systemu}) is equivalent to
\small
\begin{gather}
	\frac1{\gamma_{i^*}}\zeta_{k^*}^{i^*} e^{-\lambda \tau}u_{i^*\bar t}(k^*)-\Bigg(\frac{a_{i^*-1,k^*}-hb_{i^*k^*}}{\gamma_{i^*-1}}+\frac{-a_{i^*-1,k^*}+hb_{i^*k^*}+a_{i^*k^*}}{\gamma_{i^*+1}}\Bigg)u_{i^*x\bar x}(k^*) \nonumber \\
    -\Bigg[ (a_{i^*-1,k^*}-hb_{i^*k^*})\Big(\left(\frac1{\gamma_{i^*}}\right)_{\bar x} +\left(\frac1{\gamma_{i^*}}\right)_{x}\Big)\nonumber \\
    +\frac1{\gamma_{i^*+1}h}\Big[-a_{i^*-1,k^*}+hb_{i^*k^*}+a_{i^*k^*} \Bigg]u_{i^*\bar{x}}(k^*) \nonumber \\
    -\Bigg[(a_{i^*-1,k^*}-hb_{i^*k^*})\left(\frac1{\gamma_{i^*}}\right)_{x\bar x} +\frac1{h}(-a_{i^*-1,k^*}+hb_{i^*k^*}+a_{i^*k^*})\left(\frac1{\gamma_{i^*}}\right)_x \nonumber \\
    \frac{c_{i^*k^*}}{\gamma_{i^*}}-\frac{\zeta_{k^*}^{i^*}\lambda e^{-\lambda(t_{k^*}-t^{k^*})}}{\gamma_{i^*}}\Bigg]u_{i^*}(k^*)=e^{-\lambda t_{k^*}}f_{i^*k^*} \label{other}
\end{gather}
\normalsize

Define the sets
\begin{gather*}
\n N_+=\Bigg\{(i,k)\in\n N\Big| \\\left[ (a_{i-1,k}-hb_{ik})\Big(\left(\frac1{\gamma_i}\right)_{\bar x} + \left(\frac1{\gamma_i}\right)_{x}\Big)+\frac1{\gamma_{i+1}h}\Big[-a_{i-1,k}+hb_{ik}+a_{ik} \right] \geq 0 \Bigg\} \\ 
\n N_-=\Bigg\{(i,k)\in\n N\Big|\\ \left[ (a_{i-1,k}-hb_{ik})\Big(\left(\frac1{\gamma_i}\right)_{\bar x} + \left(\frac1{\gamma_i}\right)_{x}\Big)+\frac1{\gamma_{i+1}h}\Big[-a_{i-1,k}+hb_{ik}+a_{ik} \right] < 0 \Bigg\}.
\end{gather*}
And it's clear $\n N=\n N_+\cup\n N_-$. Suppose $(i^*,k^*)\in\n N_+$. Then owing to (\ref{systemu}) since $u_{i^*x}(k^*)\leq0$, for sufficiently small h we can write

\begin{gather}
-\Bigg[(a_{i^*-1,k^*}-hb_{i^*k^*})\Bigg(\frac1{\gamma_{i^*}}\Bigg)_{x\bar x} +\frac1{h}(-a_{i^*-1,k^*}+hb_{i^*k^*}+a_{i^*k^*})\Bigg(\frac1{\gamma_{i^*}}\Bigg)_x +\frac{c_{i^*k^*}}{\gamma_{i^*}} \nonumber\\
-\frac{\zeta_{k^*}^{i^*}\lambda e^{-\lambda(t_{k^*}-t^{k^*})}}{\gamma_{i^*}}\Bigg]u_{i^*}(k^*)\leq e^{-\lambda t_{k^*}}f_{i^*k^*}\label{est-1}
\end{gather}
If instead $(i^*,k^*)\in\n N_-$, then we can use (\ref{other}), the fact that $~u_{i^*\bar x}(k^*)\geq0$ and that
\begin{gather*}
	-\Bigg(\frac{a_{i-1,k}-hb_{ik}}{\gamma_{i-1}}+\frac{-a_{i-1,k}+hb_{ik}+a_{ik}}{\gamma_{i+1}}\Bigg)\\
	=-\frac{a_{ik}}{\gamma_{i+1}}-\frac{(\gamma_{i+1}-\gamma_{i-1})a_{i-1,k}+h(\gamma_{i-1}-\gamma_{i+1})b_{ik}}{\gamma_{i+1}\gamma_{i-1}} \\
	\leq -a_0-\frac{2h\gamma'(\tilde{x})a_{i-1,k}+h(\gamma_{i-1}-\gamma_{i+1})b_{ik}}{\gamma_{i+1}\gamma_{i-1}} \leq -\frac{a_0}{2} 
\end{gather*}
for sufficiently small h, where $\tilde{x}$ is the value that satisfies the mean value theorem to achieve again (\ref{est-1}). Therefore, (\ref{est-1}) is achieved in any case. We can choose $\tau$ so small that $e^{-\lambda(t_{k}-t^{k})}>\frac12,\quad\forall k$. The coefficient in front of $u_{i^*}(k^*)$ in \eqref{est-1} can be estimated as follows:
\begin{gather*}
	\frac{\zeta_{k^*}^{i^*}\lambda e^{-\lambda(t_{k^*}-t^{k^*})}}{\gamma_{i^*}}-(a_{i^*-1,k^*}-hb_{i^*k^*})\Bigg(\frac{-1}{\gamma_{i^*}\gamma_{i^*+1}}\gamma_{i^*x\bar{x}}+\frac{\gamma_{i^*x}+\gamma_{i^*\bar{x}}}{\gamma_{i^*-1}\gamma_{i^*}\gamma_{i^*+1}}\gamma_{i^*\bar{x}}\Bigg) \\
	+(a_{i^*k^*,\bar{x}}+b_{i^*k^*})\frac{\gamma_{i^*x}}{\gamma_{i^*}\gamma_{i^*+1}} -\frac{c_{i^*k^*}}{\gamma_{i^*}} \\
	\geq \frac{\bar{b}\lambda}{2}-(\|a\|_{L_\infty(D)}+\|b\|_{L_\infty(D)})(16\|\gamma''\|_{C[0,\ell]}+128\|\gamma'\|^2_{C[0,\ell]})\\
	-16(\|a_x\|_{L_\infty(D)}+\|b\|_{L_\infty(D)})\|\gamma'\|_{C[0,\ell]}-4\|c\|_{L_\infty(D)} \geq \frac{\bar{b}}{2}
\end{gather*}
due to definitions of $\lambda$ and $\gamma(x)$. Then by (\ref{lambda}), it is the case that the coefficient of $u_{i^*}(k^*)$ is positive independently of $i^*,k^*$. Therefore,
\[
u_{i^*}(k^*)\leq \frac{2}{\bar{b}}f_{i^*k^*}e^{-\lambda t_{k^*}}
\]

We can put together the obtained estimations to deduce that for $(i,k)\in \n M_{k_1}$,
\begin{gather*}
u_i(k)\leq\max\limits_{\n M}u_i(k) \leq A\max\Big\{0,~\Vert\Phi\Vert_{L_\infty(0,\ell)}, ~\Vert p\Vert_{L_\infty(0,T)}, \Vert g^n\Vert_{L_\infty(0,T)}, ~\Vert f\Vert_{L_{\infty}(D)}  \Big\},
\end{gather*}
with $A=\max\{1,~4a_0^{-1},~2\bar{b}^{-1}\}$.
But because $u_i(k) =\gamma_ie^{-\lambda t_k}v_i(k)$, we have the following uniform upper bound for the discrete state vector:
\begin{gather*}
v_i(k)\leq4Ae^{\lambda T}\max\Big\{0,~\Vert\Phi\Vert_{L_\infty(0,\ell)}, ~\Vert p\Vert_{L_\infty(0,T)}, ~\Vert g^n\Vert_{L_\infty(0,T)}, ~\Vert f\Vert_{L_{\infty}(D)}  \Big\},
\end{gather*}
for $\forall (i,k)\in \n M_{k_1}$. In a fully analogous manner, we arrive at a uniform lower bound for the discrete state vector: 
\begin{gather*}
v_i(k)\geq4Ae^{\lambda T}\min\Big\{0,-\Vert\Phi\Vert_{L_\infty(0,\ell)}, -\Vert p\Vert_{L_\infty(0,T)}, -\Vert g^n\Vert_{L_\infty(0,T)}, -\Vert f\Vert_{L_{\infty}(D)}  \Big\},
\end{gather*}
for $\forall (i,k)\in \n M_{k_1}$. 
Combining the uniform upper and lower bounds imply (\ref{linfestimate}) up to $k_1$. But $k_1$ was arbitrary in $1,\ldots,n$. Theorem is proved. \hfill{$\square$}

\subsection{$W_2^{1,1}(D)$ estimate for the discrete multiphase free boundary problem}

\begin{theorem}\label{energyestimate} For $[g]_n\in \n G_R^n$ and $n,m$ large enough, the discrete state vector $[v([g]_n)]_n$ satisfies the following estimate:
	\begin{gather}
	\label{energy}\Vert[v]_n\Vert^2_{\n E}:= \sum\limits_{k=1}^n\tau\sum\limits_{i=0}^{m-1}hv_{i\bar t}^2(k) +\max\limits_{1\leq k\leq n}\left(\sum\limits_{i=0}^{m-1}hv_{ix}^2(k)\right) + \sum\limits_{k=1}^n\tau^2 \sum\limits_{i=0}^{m-1} hv_{ix\bar t}^2(k)\\ \leq~~\tilde C_{\infty}\Big(\Vert\Phi\Vert^2_{W_2^1(0,\ell)} +\Vert f\Vert_{L_{\infty}(D)}^2 +\Vert p\Vert_{W_2^1(0,T)}^2 +\Vert g^n\Vert_{W_2^1(0,T)}^2\Big) \nonumber
	\end{gather}
	where $\tilde C_{\infty}$ is a constant independent of $n$ and $m$.
\end{theorem}

\noindent\emph{Proof. }Consider $n$ and $m$ large enough that Theorem \ref{uniformboundedness} is satisfied. In (\ref{dsvsum}), choose $\eta = 2\tau v_{\bar t}(k)$. Using (\ref{zetav}), write $(b_n(v_i(k)))_{\bar t} =\zeta_k^iv_{i\bar t}(k)$. Also, use the fact that
\begin{gather*}
2\tau a_{ik}v_{ix}(k)(v_{i\bar t}(k))_x \\
=a_{ik}v_{ix}^2(k)-a_{i,k-1}v_{ix}^2(k-1)-\tau a_{ik\bar{t}}v_{ix}^2(k-1)+\tau^2 a_{ik}v_{ix\bar{t}}^2(k) 
\end{gather*}
Using the above equality, and the lower bound for $a(x,t)$, we thus have
\begin{gather}
2\tau\sum\limits_{i=0}^{m-1}h\zeta_k^iv_{i\bar t}^2(k) +\sum\limits_{i=0}^{m-1}ha_{ik}v_{ix}^2(k)-\sum\limits_{i=0}^{m-1}ha_{i,k-1}v_{ix}^2(k-1)+ a_0\tau^2\sum\limits_{i=0}^{m-1}hv_{ix\bar t}^2(k)\nonumber \\ 
\leq \tau\sum\limits_{i=0}^{m-1}h a_{ik\bar{t}}v_{ix}^2(k-1)-2\tau\sum\limits_{i=0}^{m-1}hb_{ik}v_i(k)v_{ix\bar{t}}(k)-2\tau\sum\limits_{i=0}^{m-1}hc_{ik}v_i(k)v_{i\bar{t}}(k)\nonumber \\
 +2\tau\sum\limits_{i=0}^{m-1}hf_{ik}v_{i\bar t}(k)+2\tau p_kv_{m\bar t}(k)-2\tau g_k^nv_{0\bar t}(k) \label{e1}.
\end{gather}
We will now look to estimate the three summation terms on the right hand side of \eqref{e1}. By using summation by parts and Cauchy inequality with $\epsilon>0$ we get:
\begin{gather*}
	-2\tau\sum\limits_{i=0}^{m-1}hb_{ik}v_i(k)v_{ix\bar{t}}(k) = 2\tau\sum\limits_{i=1}^{m-1}hb_{ik}v_{i-1,x}(k)v_{i\bar{t}}(k)\\
	+2\tau\sum\limits_{i=i}^{m-1}hb_{i-1,k,x}v_{i-1}(k)v_{i\bar{t}}(k) - 2\tau b_{m-1,k}v_{m-1}(k)v_{m\bar{t}}(k) +2\tau b_{0k}v_0(k)v_{0\bar{t}}(k) \\
    \leq \frac{\bar{b}}{4}\tau \sum\limits_{i=0}^{m-1}hv_{i\bar{t}}^2(k)+\frac{4\|b\|^2_{L_\infty(D)}}{\bar{b}}\tau\sum\limits_{i=0}^{m-1}hv_{ix}^2(k) + \frac{\bar{b}}{4}\tau \sum\limits_{i=0}^{m-1}hv_{i\bar{t}}^2(k)\\
    +\frac{4\|b_x\|^2_{L_\infty(D)}}{\bar{b}}\tau\sum\limits_{i=0}^{m-1}hv_{i}^2(k)- 2\tau b_{m-1,k}v_{m-1}(k)v_{m,\bar{t}}(k) +2\tau b_{0k}v_0(k)v_{0\bar{t}}(k)
\end{gather*}
We will also use the fact that 
\begin{gather*}
	2\tau b_{0k}v_{0}(k)v_{0\bar{t}}(k) = \tau^2b_{0k}v_{0\bar{t}}^2(k)-b_{0k}v_{0}^2(k-1)+b_{0k}v_{0}^2(k),\\ \vspace{0.1cm} \\
    -2\tau b_{m-1,k}v_{m-1}(k)v_{m\bar{t}}(k) \\
    = -2\tau b_{m-1,k}v_{m}(k)v_{m\bar{t}}(k)+2\tau h b_{m-1,k}v_{m-1,x}(k)v_{m\bar{t}}(k) \\
    = -\tau^2b_{m-1,k}v^2_{m\bar{t}}(k)+b_{m-1,k}v^2_{m}(k-1)-b_{m-1,k}v^2_m(k)\\
    +2\tau h b_{m-1,k}v_{m-1,x}(k)v_{m\bar{t}}(k)
\end{gather*}
Estimating the other two summation terms on the right-hand side of (\ref{e1}) via Cauchy Inequality with $\ep >0$ and by recalling \eqref{b_n bound}, we have:
\begin{gather}
	2\tau\bar{b}\sum\limits_{i=0}^{m-1}hv_{i\bar t}^2(k) +\sum\limits_{i=0}^{m-1}ha_{ik}v_{ix}^2(k)-\sum\limits_{i=0}^{m-1}ha_{i,k-1}v_{ix}^2(k-1)+ a_0\tau^2\sum\limits_{i=0}^{m-1}hv_{ix\bar t}^2(k)\nonumber \\ 
    \leq \tau\sum\limits_{i=0}^{m-1}h a_{ik\bar{t}}v_{ix}^2(k-1)+ \frac{\bar{b}}{4}\tau \sum\limits_{i=0}^{m-1}hv_{i\bar{t}}^2(k)+\frac{4\|b\|^2_{L_\infty(D)}}{\bar{b}}\tau\sum\limits_{i=0}^{m-1}hv_{ix}^2(k) \nonumber\\
    + \frac{\bar{b}}{4}\tau \sum\limits_{i=0}^{m-1}hv_{i\bar{t}}^2(k) +\frac{4\|b_x\|^2_{L_\infty(D)}}{\bar{b}}\tau\sum\limits_{i=0}^{m-1}hv_{i}^2(k)+\frac{\tau}{2}\bar{b}\sum\limits_{i=0}^{m-1}hv^2_{i\bar{t}}(k)+\frac{2}{\bar{b}}\tau\sum\limits_{i=0}^{m-1}hc_{ik}^2v^2_i(k) \nonumber \\
    + \frac{\tau}{2}\bar{b}\sum\limits_{i=0}^{m-1}hv_{i\bar{t}}^2(k)+\frac{2}{\bar{b}}\tau \sum\limits_{i=0}^{m-1}hf_{ik}^2 +2\tau p_kv_{m\bar t}(k)-2\tau g_k^nv_{0\bar t}(k)+\tau^2b_{0k}v_{0\bar{t}}^2(k)\nonumber \\
    -b_{0k}v_{0}^2(k-1)+b_{0k}v_{0}^2(k) -\tau^2b_{m-1,k}v^2_{m\bar{t}}(k)+b_{m-1,k}v^2_{m}(k-1)-b_{m-1,k}v^2_m(k)\nonumber \\
    +2\tau h b_{m-1,k}v_{m-1,x}(k)v_{m\bar{t}}(k)	\label{e2}.
\end{gather}
By absorbing several terms on the right hand side of (\ref{e2}) to the left-hand side, and further bounding the right hand side we get 
\begin{gather}
	\frac{\bar{b}}{2}\tau\sum\limits_{i=0}^{m-1}hv_{i\bar t}^2(k) +\sum\limits_{i=0}^{m-1}ha_{ik}v_{ix}^2(k)-\sum\limits_{i=0}^{m-1}ha_{i,k-1}v_{ix}^2(k-1)+ a_0\tau^2\sum\limits_{i=0}^{m-1}hv_{ix\bar t}^2(k)\nonumber \\
    \leq \tau\sum\limits_{i=0}^{m-1}h a_{ik\bar{t}}v_{ix}^2(k-1)+ \frac{4\|b_x\|^2_{L_{\infty}(D)}+2\|c\|^2_{L_{\infty}(D)}}{\bar{b}}\tau\sum\limits_{i=0}^{m-1}hv^2_i(k)+\frac{2}{\bar{b}}\tau\sum\limits_{i=0}^{m-1}hf^2_{ik}\nonumber \\
    +\frac{4\|b\|^2_{L_\infty(D)}}{\bar{b}}\tau\sum\limits_{i=0}^{m-1}hv_{ix}^2(k) +2\tau p_kv_{m\bar t}(k)-2\tau g_k^nv_{0\bar t}(k)+\tau^2b_{0k}v_{0\bar{t}}^2(k)\nonumber \\
    -b_{0k}v_{0}^2(k-1)+b_{0k}v_{0}^2(k)-\tau^2b_{m-1,k}v^2_{m\bar{t}}(k)+b_{m-1,k}v^2_{m}(k-1)-b_{m-1,k}v^2_m(k)\nonumber \\
    +2\tau h b_{m-1,k}v_{m-1,x}(k)v_{m\bar{t}}(k), \label{e3}.
\end{gather}
$\forall k =\overline{1,n}$. Perform summation of (\ref{e3}) for $k$ from $1$ to $q,~~2\leq q\leq n$. The second and third terms on the left-hand side telescope, and we obtain:
\begin{gather}
	\frac{\bar{b}}{2}\sum\limits_{k=1}^{q}\tau\sum\limits_{i=0}^{m-1}hv_{i\bar t}^2(k) +\sum\limits_{i=0}^{m-1}ha_{iq}v_{ix}^2(q)+ a_0\sum\limits_{k=1}^{q}\tau^2\sum\limits_{i=0}^{m-1}hv_{ix\bar t}^2(k)\nonumber \\
    \leq \sum\limits_{k=1}^{q}\tau\sum\limits_{i=0}^{m-1}h a_{ik\bar{t}}v_{ix}^2(k-1)+\sum\limits_{i=0}^{m-1}ha_{i0}v^2_{ix}(0)\nonumber \\
    +\frac{4\|b_x\|^2_{L_{\infty}(D)}+2\|c\|^2_{L_{\infty}(D)}}{\bar{b}}\sum\limits_{k=1}^{q}\tau\sum\limits_{i=0}^{m-1}hv^2_i(k) +\frac{2}{\bar{b}}\sum\limits_{k=1}^{q}\tau\sum\limits_{i=0}^{m-1}hf^2_{ik}\nonumber \\
    +\frac{4\|b\|^2_{L_\infty(D)}}{\bar{b}}\sum\limits_{k=1}^{q}\tau\sum\limits_{i=0}^{m-1}hv_{ix}^2(k) +2\sum\limits_{k=1}^{q}\tau p_kv_{m\bar t}(k)-2\sum\limits_{k=1}^{q}\tau g_k^nv_{0\bar t}(k) \nonumber \\
	+\sum\limits_{k=1}^{q}\tau^2b_{0k}v_{0\bar{t}}^2(k) - \sum\limits_{k=1}^{q}\tau^2b_{m-1,k}v^2_{m\bar{t}}(k) + \sum\limits_{k=1}^{q} \Bigg(b_{0k}v_{0}^2(k)-b_{0k}v_{0}^2(k-1)\Bigg) \nonumber \\
    -\sum\limits_{k=1}^{q} \Bigg(b_{m-1,k}v^2_m(k)-b_{m-1,k}v^2_{m}(k-1)\Bigg)\nonumber \\
     +\sum\limits_{k=1}^{q}2\tau h b_{m-1,k}v_{m-1,x}(k)v_{m\bar{t}}(k) \label{aftertelescope}
\end{gather}
We can estimate the right hand side further and use \eqref{htau} to get
\begin{gather*}
	\sum\limits_{k=1}^{q}\tau^2b_{0k}v_{0\bar{t}}^2(k) \leq \tau \|b\|_{L_{\infty}(D)}\sum\limits_{k=1}^{q}\tau v_{0\bar{t}}^2(k) \leq \frac{\bar{b}}{8}h\sum\limits_{k=1}^{q}\tau v_{0\bar{t}}^2(k) \leq \frac{\bar{b}}{8}\sum\limits_{k=1}^{q}\tau\sum\limits_{i=0}^{m-1}hv_{i\bar t}^2(k)
\end{gather*}
Similarly,
\begin{gather*}
	- \sum\limits_{k=1}^{q}\tau^2b_{m-1,k}v^2_{m\bar{t}}(k) \leq \frac{\bar{b}}{8}\sum\limits_{k=1}^{q}\tau\sum\limits_{i=0}^{m-1}hv_{i\bar t}^2(k)
\end{gather*}
We also have:
\begin{equation}\label{est-b-esssup}
\sum\limits_{k=1}^{q} (b_{0k}v_{0}^2(k)-b_{0k}v_{0}^2(k-1)) 
	= b_{0q}v_0^2(q) - b_{00}v_0^2(0)-\sum\limits_{k=1}^{q}\tau b_{0k\bar{t}}v_0^2(k-1).
\end{equation}
Since
\begin{gather*}
\Big |\sum\limits_{k=1}^{q}\tau b_{0k\bar{t}}v_0^2(k-1)\Big |
		\leq \|[v]_n\|^2_{\ell_{\infty}}\frac{1}{h\tau} \int\limits_{0}^{h}\int\limits_{0}^{t_q}\int\limits_{t-\tau}^{t} |b_t(x,\theta)|\,d\theta \,dt\,dx\nonumber\\
		\leq \|[v]_n\|^2_{\ell_{\infty}} \int\limits_{-\tau}^{T}\esssup\limits_{x\in[0,\ell]}|b_t(x,t)|\,dt
\end{gather*} 
from \eqref{est-b-esssup} it follows that
\begin{equation*}
\Big |\sum\limits_{k=1}^{q} (b_{0k}v_{0}^2(k)-b_{0k}v_{0}^2(k-1))\Big |\leq \|[v]_n\|^2_{\ell_{\infty}}\Big( 2\|b\|_{L_{\infty}(D)}+
	\int\limits_{-\tau}^{T}\esssup\limits_{x\in[0,\ell]}|b_t(x,t)|\,dt\Big )
\end{equation*}
Similarly, we have:
\begin{gather*}
    \Big |\sum\limits_{k=1}^{q} (b_{m-1,k}v^2_m(k)-b_{m-1,k}v^2_{m}(k-1))\Big| \\
    \leq \|[v]_n\|^2_{\ell_{\infty}}\Big( 2\|b\|_{L_{\infty}(D)}+
	\int\limits_{-\tau}^{T}\esssup\limits_{x\in[0,\ell]}|b_t(x,t)|\,dt\Big )
\end{gather*}
From \cite{Abdulla5}, in a similar fashion to above, we have:
\begin{gather*}
	\sum\limits_{k=1}^{q}\tau\sum\limits_{i=0}^{m-1}h a_{ik\bar{t}}v_{ix}^2(k-1) \leq 2 \int\limits_{0}^{T}\esssup\limits_{x\in[0,\ell]}|a_t(x,t)|\,dt \max\limits_{1\leq k \leq n}\sum\limits_{i=0}^{m-1}hv_{ix}^2(k)+C\sum\limits_{i=0}^{m-1}h\Phi_{ix}^2
\end{gather*}
where $C$ is a constant independent of $n$. Using Cauchy Inequality with $\ep = \frac{\bar{b}}{8}$, we have:
\begin{gather*}
	\sum\limits_{k=1}^{q}2\tau h b_{m-1,k}v_{m-1,x}(k)v_{m\bar{t}}(k) \nonumber\\
	\leq \frac{\bar{b}}{8}\sum\limits_{k=1}^{q}\tau \sum\limits_{i=0}^{m-1} h v^2_{i\bar{t}}(k) + \frac{8\|b\|^2_{L_\infty(D)}t_q}{\bar{b}} \max\limits_{1\leq k\leq q}\sum\limits_{i=0}^{m-1} hv^2_{ix}(k)
\end{gather*}
Use the summation by parts technique on the $p$ and $g$ sums:
\begin{gather}
\sum\limits_{k=1}^q\tau p_kv_{m\bar t}(k) = 
-\sum\limits_{k=1}^{q-1}\tau p_{kt}v_m(k) + p_qv_m(q) - p_1v_m(0), \nonumber \\ 
\sum\limits_{k=1}^q\tau g_k^nv_{0\bar t}(k) = -\sum\limits_{k=1}^{q-1}\tau g_{kt}^nv_0(k) + g_q^nv_0(q) - g_1^nv_0(0) \label{sumparts}.
\end{gather}
In view of (\ref{sumparts}) and the above estimates, (\ref{aftertelescope}) yields 
\begin{gather}
	\frac{\bar{b}}{8}\sum\limits_{k=1}^{q}\tau\sum\limits_{i=0}^{m-1}hv_{i\bar t}^2(k) +\sum\limits_{i=0}^{m-1}ha_{iq}v_{ix}^2(q)+ a_0\sum\limits_{k=1}^{q}\tau^2\sum\limits_{i=0}^{m-1}hv_{ix\bar t}^2(k)\nonumber \\
	\leq 2 \int\limits_{0}^{T}\esssup\limits_{x\in[0,\ell]}|a_t(x,t)|\,dt \max\limits_{1\leq k \leq n}\sum\limits_{i=0}^{m-1}hv_{ix}^2(k)+(C+\|a\|_{L_\infty(D)})\sum\limits_{i=0}^{m-1}h\Phi_{ix}^2 \nonumber \\
    +\frac{(4\|b_x\|^2_{L_{\infty}(D)}+2\|c\|^2_{L_{\infty}(D)})T\ell}{\bar{b}}\|[v]_n\|^2_{\ell_\infty} +\frac{2}{\bar{b}}\sum\limits_{k=1}^{n}\tau\sum\limits_{i=0}^{m-1}hf^2_{ik}\nonumber\\
    \frac{4\|b\|^2_{L_\infty(D)}t_q}{\bar{b}}\max\limits_{1\leq k \leq n}\sum\limits_{i=0}^{m-1}hv_{ix}^2(k) -2\sum\limits_{k=1}^{q-1}\tau p_{kt}v_m(k) + 2p_qv_m(q) - 2p_1v_m(0)\nonumber \\
    +2\sum\limits_{k=1}^{q-1}\tau g_{kt}^nv_0(k) - 2g_q^nv_0(q) + 2g_1^nv_0(0) + \frac{8\|b\|^2_{L_\infty(D)}t_q}{\bar{b}} \max\limits_{1\leq k\leq n}\sum\limits_{i=0}^{m-1} hv^2_{ix}(k)\nonumber\\ +\Big (4\|b\|_{L_{\infty}(D)}
    + 2 \int\limits_{-\tau}^{T}\esssup\limits_{x\in[0,\ell]}|b_t(x,t)|\,dt \Big )\|[v]_n\|^2_{\ell_{\infty}}   \label{e4}
\end{gather}
Applying Cauchy-Schwartz inequality to Steklov averages, by emloying Morrey's inequality (\cite{LSU}), for sufficiently small h we have the following estimations:
\begin{gather}
\sum\limits_{i=0}^{m-1}h\Phi_{ix}^2 \leq \Vert\Phi'\Vert^2_{L_2(0,\ell)} + \Vert\Phi'\Vert_{L_2(\ell-h,\ell)}^2 \leq \Vert\Phi\Vert^2_{W_2^1(0,\ell)}, \nonumber \\
\Vert g^n\Vert_{L_{\infty}(0,T)}\leq C\Vert g^n\Vert_{W_2^1(0,T)}, \nonumber \\
\sum\limits_{k=1}^q\tau\sum\limits_{i=0}^{m-1}hf_{ik}^2 \leq \Vert f\Vert_{L_2(D)}^2, \nonumber\\
2\sum\limits_{k=1}^{q-1}\tau g_{kt}^nv_0(k) \leq \sum\limits_{k=1}^{q-1}\tau \big(g_{kt}^n\big)^2 +  \sum\limits_{k=1}^{q-1}\tau v_0^2(k)\leq  \Vert (g^n)'\Vert_{L_2(0,T)}^2+T\|[v]_n\|^2_{\ell_{\infty}} , \nonumber \\ 
2p_qv_m(q)\leq \Vert p\Vert_{L_{\infty}(0,T)}^2+v_m^2(q)\leq C\Vert p\Vert_{W_2^1(0,T)}^2+\|[v]_n\|^2_{\ell_{\infty}}, \label{estimations}.
\end{gather}
with last two estimations being extended to similar terms 
\[ 2\sum\limits_{k=1}^{q-1}\tau p_{kt}v_m(k), \ 2p_1v_m(0), \  2g_q^nv_0(q), \ 2g_1^nv_0(0)\]
Applying the results in (\ref{estimations}), along with $L_\infty$-estimate \eqref{linfestimate} to (\ref{e4}), and taking into account that  $q=\overline{1,n}$ is arbitrary we derive
\begin{gather}
	\frac{\bar{b}}{8}\sum\limits_{k=1}^n\tau\sum\limits_{i=0}^{m-1}hv_{i\bar t}^2(k) +a_0\max\limits_{1\leq k \leq n} \sum\limits_{i=0}^{m-1}hv_{ix}^2(k)+ a_0\sum\limits_{k=1}^n\tau^2 \sum\limits_{i=0}^{m-1} hv_{ix\bar t}^2(k) \nonumber \\
	\leq \tilde C_{\infty}\Big(\Vert\Phi\Vert^2_{W_2^1(0,\ell)} +\Vert f\Vert_{L_\infty(D)}^2 +\Vert p\Vert_{W_2^1(0,T)}^2 +\Vert g^n\Vert_{W_2^1(0,T)}^2\Big) \nonumber \\
	+ \Bigg(2 \int\limits_{0}^{T}\esssup\limits_{x\in[0,\ell]}|a_t(x,t)|\,dt+\frac{12\|b\|^2_{L_\infty(D)}T}{\bar{b}}\Bigg) \max\limits_{1\leq k \leq n}\sum\limits_{i=0}^{m-1}hv_{ix}^2(k) \label{e6}
\end{gather}
where $\tilde C_{\infty}$ is a constant independent of $n,m$. If 
\begin{equation}\label{smallness}
2 \int\limits_{0}^{T}\esssup\limits_{x\in[0,\ell]}|a_t(x,t)|\,dt + \frac{12\|b\|^2_{L_\infty(D)}T}{\bar{b}}< a_0
\end{equation}
then we absorb the extra term to the left-hand side and \eqref{energy} follows. If not, we partition $[0,T]$ into finitely many intervals such that in each interval $I$, \eqref{smallness} is
satisfied with integral along $I$, and $T$ in the second term replaced with the interval length $|I|$. Therefore, energy estimate \eqref{energy} holds in each interval segment, and by adding finitely many inequalities, \eqref{energy} in the whole segment follows.
\hfill{$\square$}

\subsection{Existence of Optimal Control and Convergence of Discrete Control Problems}

\begin{theorem}\label{compactness} 
    Let $\{[g]_n\}$ be a sequence in $\n G_R^n$ such that the sequence of interpolations $\{\n P_n([g]_n)\}$ converges weakly to $g\in W_2^1[0,T]$. Then the whole sequence of interpolations $\{\hat v^{\tau}\}$ of the associated discrete state vectors converges weakly in $W_2^{1,1}(D)$ to the unique weak solution $v=v(x,t;g)$ of the multiphase parabolic free boundary problem (\ref{bvpde})-(\ref{vp}).
\end{theorem}

\noindent\emph{Proof. } 
Having estimates \eqref{uniformboundedness},\eqref{energyestimate}, the proof is pursued similar to the proof of Theorem 5 in \cite{Abdulla5}. By the definitions of the interpolations in \eqref{interpolations} we have that there is a subsequence of $\{\hat{v}^\tau\}$ that converges weakly in $W_2^{1,1}(D)$ to some function $v \in W_2^{1,1}(D) \cap L_\infty(D)$, strongly in $L_2(D)$, and a further subsequence that converges to $v$ pointwise almost everywhere in $D$. We also have equivalence of $\{v^\tau\}$ and $\{\hat{v}^\tau\}$ in $W_2^{1,0}(D)$, and equivalence of $\{v^\tau\}$ and $\{\tilde{v}\}$ in $L_2(D)$. Accordingly, $v^{\tau}\rightarrow v$ weakly in $W_2^{1,0}(D)$ and $\tilde v\rightarrow v$ strongly in $L_2(D)$ and pointwise a.e. on $D$ along a subsequence. Fix arbitrary $\psi\in W_2^{1,1}(D)$ with $\psi|_{t=T}=0$. Due to density of $C^1(\bar{D})$ in $W_2^{1,1}(D)$, without loss of generality we can consider $\psi\in C^1(\overline D)$ and $\psi|_{t=T}=0$. Define $\psi_i(k) = \psi(x_i,t_k)$, and consider the interpolations:
\begin{gather}
\psi^{\tau}(x,t):= \psi_i(k),\qquad \psi_x^{\tau}(x,t):= \psi_{ix}(k)\qquad \psi_t^{\tau}(x,t):= \psi_{it}(k), \nonumber \\   x_i\leq x<x_{i+1},\quad t_{k-1}< t\leq t_{k},\qquad i =\overline{0,m},~k =\overline{0,n} \label{psi}.
\end{gather}
It is readily checked that $\psi^{\tau},\psi_x^{\tau},\psi_t^{\tau}$ converge uniformly on $\overline D$ as $n,m\rightarrow\infty$ to the functions $\psi, \psi_x,\psi_t$ respectively. For each $k$ in (\ref{dsvsum}) as satisfied by the discrete state vector $[v([g]_n)]_n$, choose $\eta_i=\tau\psi_i(k),~i=0,...,m$ and sum all equalities (\ref{dsvsum}) over $k=1,\ldots,n$. The resulting expression is as follows:
\begin{gather}
    \sum\limits_{k=1}^n\tau\sum\limits_{i=0}^{m-1}h\Big[\big(b_n(v_i(k))\big)_{\bar t}\psi_i(k) + a_{ik}v_{ix}(k)\psi_{ix}(k)+b_{ik}v_i(k)\psi_{ix}(k)\nonumber \\
    +c_{ik}v_i(k)\psi_i(k)-f_{ik}\psi_i(k)\Big]-\sum\limits_{k=1}^n\tau p_k\psi_m(k) +\sum\limits_{k=1}^n\tau g_k^n\psi_0(k) = 0 \label{c1}
\end{gather}
The first term is transformed through summation by parts as in \cite{Abdulla5}, and using the interpolations, \eqref{c1} becomes the following integral identity:
\begin{gather}
    \int\limits_0^T\int\limits_0^{\ell}\Big[-b_n(\tilde v)\psi_t^{\tau} +av_x^{\tau}\psi_x^{\tau}+b\tilde{v}\psi_x^{\tau}+c\tilde{v}\psi^\tau -f\psi^{\tau}\Big]\,dx\,dt -\int\limits_0^{\ell}b_n(\tilde\Phi)\psi^{\tau}(x,\tau)\,dx  \nonumber \\ -\int\limits_0^Tp(t)\psi^{\tau}(\ell,t)\,dt +\int\limits_0^Tg^n(t)\psi^{\tau}(0,t)\,dt + \int\limits_{T-\tau}^T\int\limits_0^{\ell}b_n(\tilde v)\psi_t^{\tau}\,dx\,dt = 0 \label{c2}.
\end{gather}
From \cite{Abdulla5}, we have that $b_n(\tilde{v})$ has a subsequence such that both $b_n(\tilde{v})$ and $b_n(\tilde{\Phi})$ converge weakly in $L_2(D)$ and $L_2[0,\ell]$, respectfully, to functions of type $\n B$, which we will denote as $\tilde{b}(x,t)$ and $\tilde{b}_0(x)$. We also have that the last integral tends to $0$ due to absolute continuity of the integral. Using the convergence properties of the interpolations, 
due to weak convergence of $\{\hat{v}^\tau\}$, equivalence of $\{\hat{v}^\tau\}$ and $\{v^\tau\}$, and uniform convergence of $\{\psi^\tau_x\}$, passing to the limit as $n\to+\infty$ we get:
\begin{gather}
    \int\limits_0^T\int\limits_0^{\ell}\Big[-\tilde b(x,t)\psi_t +av_x\psi_x+bv\psi_x+cv\psi -f\psi\Big]\,dx\,dt -\int\limits_0^{\ell}\tilde b_0(x)\psi(x,0)\,dx  \nonumber \\ -\int\limits_0^Tp(t)\psi(\ell,t)\,dt +\int\limits_0^Tg(t)\psi(0,t)\,dt = 0 \label{vlimweaksol}.
\end{gather}
Since $\tilde{b}(x,t)$ and $\tilde{b}_0(x)$ are both of type $\n B$, and by use of Mazur's lemma, we deduce as in \cite{Abdulla5} that $\tilde{b}(x,t) = B(x,t,v(x,t))$ $\tilde{b}_0(x)=B(x,0,\Phi(x))$ a.e on $D$ and $(0,\ell)$ respectfully. This implies that $v$ is a weak solution in the sense of Definition~\ref{weaksol}. By Lemma \ref{unique}, this implies $v$ is the only solution of the problem, and hence the only limit point of the sequence $\{\hat{v}^\tau\}$. \hfill{$\square$}

Having estimates  \eqref{uniformboundedness},\eqref{energyestimate} and compactness Theorem~\ref{compactness}, the completion of the proof of Theorems~\ref{existence} and~\ref{convergence}  coincides with the proof given in \cite{Abdulla5}, through compactness arguments and proving weak continuity of cost functional $\n J$, and verification of the conditions of Theorem~\ref{Vasil} (see Lemmas A, B and C in \cite{Abdulla5}).

\bibliographystyle{amsplain}

\end{document}